\input amstex
\documentstyle{amsppt}
\magnification=\magstep1 \NoRunningHeads
\topmatter

\title
K-property for~Maharam~extensions~of  nonsingular Bernoulli and  Markov shifts
 \endtitle

\author
Alexandre I. Danilenko and Mariusz Lema{\'n}czyk
\endauthor

\thanks{Research of the second named author is supported in part by Narodowe Centrum Nauki  UMO-2014/15/B/ST1/03736.}
\endthanks
\thanks{Research of the two authors is supported in part by the special program of invitations of the semester
``Ergodic Theory and Dynamical Systems in their Interactions with Arithmetic and Combinatorics'', Chair Jean Morlet, 1.08.2016-30.01.2017.}
\endthanks

\email
alexandre.danilenko@gmail.com
\endemail

\address
 Institute for Low Temperature Physics
\& Engineering of National Academy of Sciences of Ukraine, 47 Nauky Ave.,
 Kharkiv, 61103, UKRAINE
\endaddress
\email alexandre.danilenko\@gmail.com
\endemail

\address
Faculty of Mathematics and Computer Science, Nicolaus Copernicus
University, ul. Chopina 12/18, 87-100 Toru\'n, Poland
\endaddress

\email
mlem\@mat.umk.pl
\endemail

\abstract
It is shown that each  conservative nonsingular Bernoulli shift is either of type $II_1$ or $III_1$.
Moreover, in the latter case the corresponding Maharam extension of the shift  is a $K$-automorphism.
This extends earlier results obtained by Z.~Kosloff for the equilibrial shifts.
Nonequilibrial shifts of type $III_1$ are constructed.
We further generalize  (partly) the main results  to  nonsingular Markov shifts.
\endabstract

 \loadbold

\NoBlackBoxes

\endtopmatter

\document

\head 0. Introduction
\endhead

In this paper we study asymptotic properties of nonsingular Bernoulli and nonsingular Markov shifts.
By a {\it nonsingular Bernoulli shift} we mean the 2-sided shift $\widetilde T$ on the infinite product space $\widetilde X:=\{0,1\}^\Bbb Z$ equipped with the infinite product measure $\bigotimes_{i\le 0}\mu_1\otimes\bigotimes_{i\ge 1}\mu_i$, where the probability measures $\mu_i$, $i\ge 1$, are chosen in  such a way that the shift is non-singular.
We call the Bernoulli shift  {\it equilibrial}  if $\mu_1(0)=\mu_1(1)=0.5$.
We are interested in nonsingular Bernoulli shifts not admitting finite invariant equivalent measures.
It is easy to construct dissipative (in fact, totally dissipative) nonsingular Bernoulli shifts.
The first example of a conservative nonsingular Bernoulli shift not admitting a  finite invariant equivalent measure was constructed by Krengel  in 1970 \cite{Kr}.
Later Hamachi presented another family of conservative nonsingular Bernoulli shifts of type $III$, i.e. transformations which  have no $\sigma$-finite invariant equivalent measures, neither finite nor infinite \cite{Ha}.
This was even more refined in \cite{Ko1}, where Kosloff constructed a family of nonsingular Bernoulli shifts of type $III_1$.
This means that the Maharam extension of the shift is ergodic.
In a subsequent paper \cite{Ko2} he showed that each conservative nonsingular  equilibrial Bernoulli shift either admits a finite  invariant equivalent measure or is of type $III_1$. Moreover, in the latter case the corresponding Maharam extension of the shift has property $K$ (in the sense of Silva and Thieullen \cite{SiTh}).
This implies that  the aforementioned Bernoulli shifts from \cite{Kr} and \cite{Ha} are all of type $III_1$.

Only equilibrial Bernoulli shifts were considered in \cite{Kr}, \cite{Ha}, \cite{Ko1} and \cite{Ko2}.
In the first part of the present paper we extend the main results of \cite{Ko2} to the general nonsingular Bernoulli shifts:

\proclaim{Theorem 3.1 and Corollary 3.3}
\roster
\item"$(i)$"
If $\sum_{n=1}^\infty(\mu_n(0)-\mu_1(0))^2<\infty$ then $\mu$ is equivalent to 
$\bigotimes_{n=-\infty}^\infty\mu_1$ and hence
$\widetilde T$ is of type $II_1$. 
\item"$(ii)$"
If $\sum_{n=1}^\infty(\mu_n(0)-\mu_1(0))^2=\infty$ and $\widetilde T$ is conservative  then $\widetilde T$ is ergodic of type $III_1$ and the Maharam extension  of $\widetilde T$ is a weakly mixing  $K$-automorphism.
\endroster
\endproclaim

The main difference between the general case  and the equilibrial one  is that in the general case the cocycle generated by the Radon-Nikodym derivative of a one-sided Bernoulli shift and the Radon-Nikodym cocycle of the tail equivalence relation generated by the shift do not coincide (see the proof of Theorem~3.1).
Hence we can not apply the Araki-Woods lemma on the Krieger's type of the ITPFI factors \cite{ArWo} which played a crucial role in the equilibrial case.
That is why we  need first to establish a stronger version of the Araki-Woods lemma (see Proposition~1.5).

We also give a simple  explicit inductive construction of the sequence of measures $(\mu_i)_{i=1}^\infty$ with an arbitrary $\mu_1$ such that the corresponding Bernoulli shift is nonsingular, conservative and not of type $II_1$ (see Theorem~4.2).
The conservativeness is the key property to establish.
We achieve it by approximating  the shift with a sequence of Bernoulli type $II_1$ shifts
each of which is, of course, conservative.
We do  the approximation in such a way  that $\widetilde T$ inherits conservativeness partly, i.e. on a finite subalgebra of cylinders,   from each of these prelimit transformations.
In the limit, the subalgebras generate the entire Borel $\sigma$-algebra.
Hence $\widetilde T$ is conservative.

In the second part of the paper we consider {\it nonsingular Markov shifts}, i.e. 2-sided shifts on $\widetilde X$ equipped with Markov measure $\widetilde\mu$ determined by a probability $\lambda$ on $\{0,1\}$ and a sequence of stochastic matrices $(P^{(n)})_{n=1}^\infty$ (see Sections~5 and 7). 
In the bistochastic case we prove the following analogue of Theorem~3.1.

\proclaim{Theorem 8.1} 
Let the matrices $P^{(n)}$, $n\ge 1$, be all bistochastic,
  $\lambda(0)=\lambda(1)=0.5$ and $P^{(1)}=\pmatrix 0.5 &0.5\\0.5&0.5\endpmatrix$.
If $(\widetilde X,\widetilde\mu,\widetilde T)$ is conservative then it is weakly mixing and either of type $II_1$ (if  $\sum_{n\ge 1}(P^{(n)}_{0,0}-0.5)^2<\infty$) or of type $III_1$ (otherwise).
 In the latter case,
the Maharam extension of $(\widetilde X,\widetilde{\mu},\widetilde T)$ is a $K$-automorphism.
\endproclaim

We also prove some partial analogues of Theorem~3.1 in the general (not only bistochastic) Markov case in Theorems~7.5 and 9.4.

The outline of the paper is as follows.
In Section~1 we first briefly remind the basic concepts of measurable orbit theory: nonsingular equivalence relation, skew-product extension, essential value of a cocycle, Krieger's type of an equivalence relation, etc. Then we prove some generalizations of Araki-Woods lemma (see Proposition~1.5, Remarks~1.6 and 1.7) that will be utilized in the later sections.
In Section~2 we review the general theory of nonsingular endomorphisms and its relation to the measurable orbit theory.
We collect there some  facts about conservativeness, recurrence, ergodicity,  exactness, Maharam extension and natural extension for endomorphisms.
In Section~3 we prove one of the main results of the paper.
The Maharam extension of a conservative nonsingular  Bernoulli shift (which is the natural extension of a one-sided nonsingular Bernoulli shift admitting no equivalent invariant probability measure) is a $K$-automorphism (Theorem~3.1).
In Section~4 we give  concrete examples of  conservative nonsingular Bernoulli shifts of type $III_1$.
In Section~5 we introduce Markov measures on $\{0,1\}^\Bbb Z$.
Some analogues of Kakutani theorem on equivalence of infinite product measures 
and Kolmogorov zero-one law for the Markov measures
are under discussion there.
In Section~6 we compute Krieger's type of the tail equivalence relation equipped with a stationary Markov measure.
In Section~7 we introduce nonsingular one-sided Markov shifts and describe their natural extensions. We  find a sufficient condition for the natural extensions of the shifts to be  $K$-automorphisms (Theorem~7.5).  A necessary condition for conservativeness of the natural extensions  is also found (Lemma~7.6).
In Section~8 we consider bistochastic nonsingular Markov shifts.
We show how to extend the main results obtained  in Section~3 for  Bernoulli shifts to the bistochastic Markov shifts under some ``initial conditions"  (Theorem~8.1).
The general Markov case is considered in Section~9 (see Theorem~9.4).
Section~10 is a list of open problems and remarks related to the subject of the paper.

\head 1. Measurable equivalence relations and~their~cocycles.
Generalizations of Araki-Woods lemma
\endhead

Let $(X,\goth B,\mu)$ be a standard $\sigma$-finite measure space.
A Borel equivalence relation $\Cal R\subset X\times X$ is called {\it countable} if for each $x\in X$, the $\Cal R$-equivalence class $\Cal R(x)$ is countable.
$\Cal R$ is called {\it $\mu$-nonsingular} if for each subset $A\in\goth B$ of zero measure, the subset $\Cal R(A):=\bigcup_{x\in A}\Cal R(x)$ is also of zero measure.
If for each $A\in\goth B$ of positive measure  the intersection $(\Cal R(x)\setminus \{x\})\cap A$ is nontrivial for a.e. $x\in A$ then  $\Cal R$ is called {\it $\mu$-conservative}.
If the $\sigma$-algebra of $\Cal R$-invariant (i.e. $\Cal R$-saturated) Borel subsets in $X$ is trivial (mod $\mu$) then $\Cal R$ is called {\it $\mu$-ergodic}.

From now on we will assume that $\Cal R$ is countable and $\mu$-nonsingular.
Given a locally compact second countable  group $G$, a Borel map $\alpha:\Cal R\to G$ is called a {\it cocycle} of $\Cal R$
if there is a $\mu$-conull subset $A\subset S$ such that $\alpha(x,y)\alpha(y,z)=\alpha(x,z)$ for all $x,y,z\in A$ such that $(x,y),(y,z)\in\Cal R$.
If  $\lambda_G$ is a left Haar measure on $G$ then we can define the {\it $\alpha$-skew product}  equivalence relation $\Cal R(\alpha)$ on the product space $(X\times G,\mu\times\lambda_G)$ by setting
$(x,g)\sim(y,h)$ if $(x,y)\in\Cal R$ and $h=\alpha(x,y)g$.
Of course, $\Cal R(\alpha)$ is countable and $(\mu\times\lambda_G)$-nonsingular.
If $\Cal R(\alpha)$ is conservative then  $\alpha$ is called {\it recurrent}.
If $\Cal R(\alpha)$ is ergodic then  $\alpha$ is called {\it ergodic}.
Two cocycles $\alpha,\beta:\Cal R\to G$ are called {\it cohomologous} if there is a Borel map $\phi:X\to G$ such that $\alpha(x,y)=\phi(x)^{-1}\beta(x,y)\phi(y)$ for all $(x,y)\in\Cal R\cap(B\times B)$, where $B\subset X$ is a $\mu$-conull subset.
A cocycle is called a {\it  coboundary} if it is cohomologous to the trivial cocycle.

 The {\it  Radon-Nikodym} cocycle $\Delta_{\Cal R,\mu}:\Cal R\to\Bbb R^*_+$ (of the  pair $(\Cal R,\mu)$) can be defined in several ways.
We consider  only the following one.
Let $\Gamma$ be a countable group of Borel bijections of $X$ such that $\Cal R(x)=\{\gamma x\mid\gamma\in\Gamma\}$ for each $x\in X$.
It exists (but is non-unique) according to \cite{FeMo}.
We now set
$$
\Delta_{\Cal R,\mu}(x,\gamma x):=(d\mu\circ \gamma/d\mu)(x),\qquad x\in X.
$$
 Then $\Delta_{\Cal R,\mu}$ is well defined (does not depend on the particular choice of $\Gamma$).
 If $\nu$ is a $\sigma$-finite measure on $X$ which is equivalent to $\mu$ then $\Delta_{\Cal R,\mu}$
 is cohomologous to $\Delta_{\Cal R,\nu}$.
 Conversely,  a cocycle of $\Cal R$ which is cohomologous to $\Delta_{\Cal R,\mu}$ is $\Delta_{\Cal R,\nu}$ for a measure $\nu$ equivalent to $\mu$.
By the Maharam theorem, $\Cal R$ is $\mu$-conservative if and only if $\Delta_{\Cal R,\mu}$ is recurrent \cite{Sc}.
We say that $\mu$ is {\it $\Cal R$-invariant} if $\Gamma$ preserves $\mu$ (this does not depend on the particular choice of $\Gamma$).

Suppose that $\Cal R$ is $\mu$-ergodic.
Given a cocycle $\alpha$ of $\Cal R$ with values in an Abelian group $G$, an element $g\in G$ is called an {\it essential value} of $\alpha$ if for each set $A\in\goth B$ of positive measure and each neighborhood $U$ of $g$ in $G$, there is a subset $B\subset A$ of positive measure and a one-to-one  mapping $\gamma:B\to A$ such that $(x,\gamma x)\in\Cal R$ and $\alpha(x,\gamma x)\in U$ for each $x\in B$.
The set of all essential values of $\alpha$ is denoted by $r(\alpha)$.
It is a closed subgroup of $G$.
If  a cocycle $\beta:\Cal R\to G$ is cohomologous to $\alpha$ then
$r(\alpha)=r(\beta)$.
The cocycle $\alpha$ is ergodic if and only if $r(\alpha)=G$ \cite{Sc}.
It is easy to verify that given another Abelian locally compact second countable group $H$ and a homomorphism $\theta:G\to H$, then $\theta(r(\alpha))\subset r(\theta\circ\alpha)$.

In order to verify that an element of $G$ is an essential value of  $\alpha$ we will use the following approximation lemma.

\proclaim{Lemma 1.1 \rom{(cf. \cite{Ch--Pr, Lemma~2.1})}}
Let $\goth A\subset\goth B$ be a semiring  such that the corresponding ring  $\goth F(\goth A) $  is dense in $\goth B$.
Let $1>\delta>0$ and let $g\in G$.
If for each set $A\in\goth A$ of positive measure and a neighborhood $U$ of $g$ there is a subset $B\subset A$ and a one-to-one mapping $\gamma:B\to A$ such that  $\mu(B)>\delta\mu(A)$, $(x,\gamma x)\in\Cal R$, $\alpha(x,\gamma x)\in U$ and $\delta<\Delta_{\Cal R,\mu}(x,\gamma x)<\delta^{-1}$
for all $x\in B$ then $g\in r(\alpha)$.

\endproclaim

Suppose that  $\mu$ is non-atomic and $\Cal R$ is $\mu$-ergodic.
Then $\Cal R$ is called {\it of type $II$} if there is a $\sigma$-finite $\Cal R$-invariant measure $\nu$ equivalent to $\mu$ or, equivalently, $\Delta_{\Cal R,\mu}$ is a coboundary.
If, moreover, $\nu(X)<\infty$ then $\Cal R$ is called {\it of type $II_1$}.
If $\nu(X)=\infty$ then $\Cal R$ is called {\it of type $II_\infty$}.
If $\Cal R$ is not of type $II$ then $\Cal R$ is called of type $III$.
The type $III$ admits further classification into subtypes $III_\lambda$, $0\le\lambda\le 1$.
If $\Delta_{\Cal R,\mu}$ is ergodic, i.e. $r(\Delta_{\Cal R,\mu})=\Bbb R^*_+$, then $\Cal R$ is called {\it of type $III_1$}.
If there is  $\lambda\in (0,1)$ such that $\Delta_{\Cal R,\mu}$ is cohomologous to a cocycle  $\beta$ taking  values in the closed subgroup $\{\lambda^n\mid n\in\Bbb Z\}\subset\Bbb R^*_+$ and $\beta$ is ergodic as a cocycle with values in this subgroup then
$\Cal R$ is called {\it of type $III_\lambda$}.
Finally, if $\Cal R$ is if type $III$ but not of type $III_\lambda$, $0<\lambda\le 1$, then $\Cal R$ is called {\it of type $III_0$}.
Equivalently, $\Cal R$ is of type $III_0$ if $\Cal R$ is of type $III$ and $\Delta_{\Cal R,\mu}=\{1\}$ (we refer to \cite{FeMo} and \cite{HaOs1} for details).

We  recall an easy  classic  lemma on the type of direct product of two ergodic equivalence relations.

\proclaim{Lemma 1.2 (\cite{ArWo}, \cite{HaOs1})}
Let $\Cal R_i$ be an ergodic $\mu_i$-nonsingular equivalence relations on a standard measure space $(X_i,\goth B_i,\mu_i)$, $i=1,2$.
\roster
\item"$(i)$"
If $\Cal R_1$ is of type $II_1$ then $\Cal R_1\times\Cal R_2$ is of the same type as $\Cal R_2$.
\item"$(ii)$"
If $\Cal R_1$ is of type $III_1$ then so is $\Cal R_1\times\Cal R_2$.
\item"$(iii)$"
Let $\Cal R_1$ be of type $III_\lambda$ and let $\Cal R_2$ be of type $III_\xi$
with $0<\lambda,\xi<1$ and let $\Lambda:=\{\lambda^n\xi^m\mid n,m\in\Bbb Z\}$.
 If $\Lambda=\{\eta^n\mid n\in\Bbb Z\}$ for some $\eta\in(0,1)$ then
$\Cal R_1\times\Cal R_2$ is of type $III_\eta$.
Otherwise $\Cal R_1\times\Cal R_2$ is of type $III_1$.
\endroster
\endproclaim

In what follows we will consider infinite product spaces.
Let $A$ be a finite set.
Then the space  $A^\Bbb N$ endowed with the topology of  infinite product of the discrete topologies is a compact metric space.
Given $n\le m$ and a finite sequence $a_n,\dots,a_m$ of elements from $A$,
we denote by $[a_n,\dots,a_m]_n^m$ the corresponding {\it cylinder} in $A^\Bbb N$, i.e. the subset $\{x=(x_j)_{j=1}^\infty\in A^\Bbb N\mid x_i=a_i\text{ whenever }n\le i\le m\}$.
Since the algebra  $\goth K$ consisting  of the finite unions of cylinders is nothing but the algebra of clopen subsets in $A^\Bbb N$ which is  a base of the topology
on $A^\Bbb N$, it  follows that $\goth K$ is dense in the $\sigma$-algebra
 of Borel subsets with respect to any measure on $A^\Bbb N$.
The  {\it tail equivalence relation} on $A^\Bbb N$ is defined by
$$
(x_i)_{i=1}^\infty\sim (y_i)_{i=1}^\infty \text{ if there is }N>0\text{  such that $x_i=y_i$ for all $i\ge N$.}
$$
Let $(X,\mu)=(\{0,1\}^\Bbb N,\mu_1\times\mu_2\times\cdots)$  such that  $\mu_i(0)>0$, $\mu_i(1)>0$  and $\mu_i(0)+\mu_i(1)=1$ for all $i\in\Bbb N$.
Let $\Cal R$ denote the  tail equivalence relation on $X$.
It is easy to verify that  $\Cal R$ is $\mu$-nonsingular and
$\Delta_{\Cal R,\mu}(x,y)=\prod_{i>0}\mu_i(y_i)/\mu_i(x_i)$ for all  pairs $(x,y)\in\Cal R$, $x=(x_i)_{i>0},y=(y_i)_{i> 0}$.
According to Kolmogorov zero-one law, $\Cal R$ is $\mu$-ergodic.

\proclaim{Lemma 1.3 \cite{ArWo}}
Let $(X,\mu,\Cal R)$ be as above and
 there exist  $\lambda\in (0,1]$ and a sequence $\epsilon_i\to 0$ such that
$$
\mu_i(0)=\frac 1{1+\lambda e^{\epsilon_i}}\quad\text{and}\quad  \mu_i(1)=\frac {\lambda e^{\epsilon_i}}{1+\lambda e^{\epsilon_i}}, \quad i>0.
\tag1-1
$$
\roster
\item"$(i)$"
If \,
$\sum_{i>0}\epsilon_i^2<\infty$ then
$\mu$ is equivalent to the infinite product $\kappa_\lambda\times\kappa_\lambda\times\cdots$, where $\kappa_\lambda(0)=1/(1+\lambda)$ and $\kappa_\lambda(1)=\lambda/(1+\lambda)$.
Hence
$\Cal R$ is of type $III_\lambda$ in case $0<\lambda<1$ or of type $II_1$ in case $\lambda=1$.
\item"$(ii)$"
If \,
$\sum_{i>0}\epsilon_i^2=\infty$ then $\Cal R$ is of type $III_1$. 
\endroster
\endproclaim

\remark{Remark 1.4}
The claim $(i)$ of  Lemma 1.3 follows straightforwardly  from the Kakutani theorem on equivalence of infinite product measures \cite{Ka}.
The claim $(ii)$ is more involved.
It is  a particular case of  Lemma~9.3 from the paper \cite{ArWo}
devoted  to  the theory of  operator algebras.
It was an attempt in \cite{Os} to give a pure measure theoretical proof of this result.
However, in our opinion, that proof  has a couple of  flaws (for instance, in the place where the author applies the central limit theorem).
 In  \cite{BrDo} Brown and Dooley used the language of $G$-measures
to
provide a new proof of  a simple case of Lemma~1.3$(ii)$  (with $\lambda=1$).
In a subsequent paper \cite{Br--La}  Brown, Dooley and Lake showed that
this proof is false.
They  gave another proof  (only for the case $\lambda=1$) which  did not use the $G$-measures \cite{Br--La}.
In our opinion, their new proof is somewhat more complicated  comparatively with the original  argument by Araki and Woods.
Below we will need the following proposition whose proof implies Lemma~1.3$(ii)$ (see Remark~1.6).
To prove it we use an argument which is  close to the  argument utilized in \cite{ArWo, Lemma~9.3}.
\endremark

\proclaim{Proposition 1.5}
Let $(X,\mu,\Cal R)$ be as in Lemma~1.3(ii).
Define a cocycle $\Lambda:\Cal R\to\Bbb R$  by setting
$$
\Lambda(x,y)=\sum_{i>0}(y_i-x_i)\log\lambda.
$$
Then the cocycle  $\alpha:=\log\Delta_{\Cal R,\mu}-\Lambda$ of $\Cal R$ with values in $\Bbb R$
 is ergodic.
\endproclaim

\demo{Proof}
Fix an infinite subset
$J\subset\Bbb N$ such that
$\sum_{i\in J}{\epsilon_i^2}<\infty$.
Applying Lem\-ma~1.3(i) we  replace
 $\mu$ by an equivalent measure for which \thetag{1-1} is satisfied, $\sum_{i>0}\epsilon^2_i=\infty$ and $\epsilon_i=0$ if $i\in J$.
 Therefore without loss of generality we may think that the triplet $(X,\mu,\Cal R)$ is isomorphic to the  triplet
 $
(Z, \eta,\Cal T),
 $
 where $Z:=(\{0,1\}\times\{0,1\})^\Bbb N$, $\eta:=\bigotimes_{i=1}^\infty (\mu_i\times\mu_0)$, $\mu_0(0):=1/(1+\lambda)$ and $\mu_0(1):=\lambda/(1+\lambda)$ and $\Cal T$ is the tail equivalence relation on $Z$.
 It is easy to verify that
 $$
\log \Delta_{\Cal T,\eta}(z,z')=\sum_{i>0}\left((x_i'-x_i)(\log\lambda+\epsilon_i)+(y_i'-y_i)\log\lambda
\right)
\tag1-2
 $$
  where $z=(x_i,y_i)_{i>0}$ and  $z'=(x_i',y_i')_{i>0}$ are $\Cal T$-equivalent points, and $x_i,x_i',y_i,y_i'\in\{0,1\}$ for each $i>0$.
Computing $\alpha$ in the ``new coordinates'' we obtain that
$$
\alpha(z,z')=\sum_{i>0}\epsilon_i(x_i'-x_i)\qquad\text{for all }(z,z')\in\Cal T.\tag1-3
$$
Denote by $\tau$ the flip on $\{0,1\}\times\{0,1\}$, i.e. $\tau(i,j)=(j,i)$.
Given $z=(x_i,y_i)_{i>0}\in Z$ and $n>0$, we denote by $z^{*n}$ the element $(x_i^*,y_i^*)_{i>0}\in Z$, where
$$
(x_i^*,y_i^*):=\cases
\tau(x_i,y_i), &\text{for }i=1,\dots,n\\
(x_i,y_i),&\text{otherwise}.
\endcases
$$
Of course, $(z,z^{*n})\in\Cal T$.

{\it Claim 1.} For each $a\in\Bbb R$,
$$
\lim_{n\to\infty}\eta(\{z\in Z\mid \alpha(z,z^{*n})> a\})=0.
$$
To prove this claim we define mappings $X_i:Z\to\Bbb R$ by
$$
X_i(z):=
\cases
(-1)^{x_i}\epsilon_i, &\text{if  }x_i\ne y_i\\
0,&\text{otherwise.}
\endcases
$$
Then $X_1,X_2,\dots$ is a sequence of independent random variables and $|X_i|\le 1$ for all $i$.
The expected value $E(X_i)$ equals
 $$
\epsilon_i\mu_i(0)\mu_0(1)-\epsilon_i\mu_i(1)\mu_0(0)=
\frac{\lambda\epsilon_i(1- e^{\epsilon_i})}{(1+\lambda)(1+\lambda e^{\epsilon_i})}=
\frac{-\lambda\epsilon_i^2}{(1+\lambda)^2}(1+\overline{o}(1)).
$$
In a similar way,
$$
E(X_i^2)=
\frac{\lambda\epsilon_i^2(1+ e^{\epsilon_i})}{(1+\lambda)(1+\lambda e^{\epsilon_i})}=
\frac{2\lambda\epsilon_i^2}{(1+\lambda)^2}(1+\overline{o}(1)).
$$
 Therefore the variance $\sigma^2(X_i)$ equals
$\frac{2\lambda\epsilon_i^2}{(1+\lambda)^2}(1+\overline{o}(1))$.
Hence $\lim_{i\to\infty}E(X_i)=0$ and $\sum_{i=1}^\infty\sigma^2(X_i)=+\infty$.
It now follows from the central limit theorem for uniformly bounded sequences of independent random variables that
$$
\frac{\sum_{i=1}^n(X_i-E(X_i))}{\sqrt{\sum_{i=1}^n\sigma^2(X_i)}}\quad\text{approaches the normal distribution as $n\to\infty$.}\tag 1-4
$$
We note that
$$
\alpha(z,z^{*n})=\sum_{i=1}^n\epsilon_i(y_i-x_i)=\sum_{i=1}^nX_i(z).
$$
Let $A_n:=\{z\in Z\mid \alpha(z,z^{*n})>a \}$.
Since
 $\sqrt{\sum_{i=1}^n\sigma^2(X_i)}\to\infty$
and
$$
-\infty =
\lim_{n\to\infty}\frac{-\sqrt\lambda\sum_{i=1}^n\epsilon_i^2}{(1+\lambda)\sqrt{2\sum_{i=1}^n\epsilon_i^2}}=
 \lim_{n\to\infty}\frac{\sum_{i=1}^nE(X_i)}{\sqrt{\sum_{i=1}^n\sigma^2(X_i)}},
$$
it follows that
$$
\min_{z\in A_n}\frac{\sum_{i=1}^n(X_i(z)-E(X_i))}{\sqrt{\sum_{i=1}^n\sigma^2(X_i)}}\to+\infty
$$
as $n\to\infty$.
In view of \thetag{1-4}, we obtain that
 $\eta(A_n)\to 0$, as claimed.

{\it Claim 2}.
Fix $r<-1$.
We are going to show that $r$ is an essential value of $\alpha$.
Let $\epsilon>0$.
Choose $k>0$ such that $\epsilon_i<\epsilon$ for all $i\ge k$.
Fix a cylinder $C:=[(a_1,b_1),\dots,(a_k,b_k)]_{1}^k\subset Z$.
It follows from the proof of Claim~1 that there are $N>0$ and a subset $I\subset (\{0,1\}\times\{0,1\})^{N-k}$ such that
for the subset
$$
A:=\bigsqcup_{(z_{k+1},\dots,z_N)\in I}[(a_1,b_1),\dots,(b_k,a_k),z_{k+1},\dots,z_N]_1^N\subset C,
$$
we have
$
\eta(A)>0.5\eta(C)
$
and
$\max_{z\in A}\alpha(z,z^{\bullet N})< r$,
where
$$
z^{\bullet l}:=((a_1,b_1),\dots,(a_k,b_k),\tau(z_{k+1}),\dots,\tau(z_l),z_{l+1},z_{l+2},\dots)
$$
and $l=k+1,\dots,N$.
Note that by \thetag{1-3}, $| \alpha(z,z^{\bullet l})-\alpha(z,z^{\bullet (l+1)})|\le\epsilon_{l+1}$.
For $z\in A$, let $l(z)$ be the smallest number $l>k$ such that
$\alpha(z,z^{\bullet l})<r$.
Then $l(z)\le N$ and
$
|\alpha(z,z^{\bullet l(z)})-r|\le \epsilon_{l(z)}<\epsilon.
$
We now set
$$
\phi(z):=z^{\bullet l(z)},\quad z\in A.
$$
Then  $(z,\phi(z))\in\Cal T$, the mapping $\phi:A\ni z\mapsto\phi(z)\in C$ is one-to-one and
$|\alpha(z,\phi(z))-r|<\epsilon$ for all $z\in A$.
To show that $\phi$ is one-to-one, we suppose that $\phi(z)=\phi(z')$ for some
$z=(z_i)_{i=1}^\infty,z'=(z'_i)_{i=1}^\infty\in A$.
If $l(z)=l(z')$ then obviously $z=z'$.
Therefore suppose that $l(z)>l(z')$.
Then the equality $\phi(z)=\phi(z')$ yields that  $z_i=z_i'$ if $1\le i\le k$ and
$\tau(z_i)=\tau(z_i')$ if  $k+1\le i\le l(z')$.
Hence $z_i=z_i'$  whenever $1\le i\le  l(z')$.
Therefore $\alpha(z,z^{\bullet l(z')})=\alpha(z',(z')^{\bullet l(z')})<r$, a contradiction which implies $l(z)=l(z')$ and hence $z=z'$.
It follows from \thetag{1-2}  and \thetag{1-3} that $\log\Delta_{\Cal T,\eta}(z,\phi(z))=\alpha(z,\phi(z))$ for all $z\in A$.
Now  Lemma~1.1 yields (put $\delta:=e^{r-\epsilon}$) that $r$ is an essential value of $\alpha$.
Since $r$ is an arbitrary  real number  less than $-1$, it follows that  $r(\alpha)=\Bbb R$, i.e. $\alpha$ is ergodic.
\qed
\enddemo

\remark{Remark 1.6} In fact, we proved  more than  claimed in the statement of  Proposition~1.5.
It was shown indeed that $(r,r)$ is an essential value for the ``double cocycle" $\alpha\times\log\Delta_{\Cal R,\mu}:\Cal R\to\Bbb R\times\Bbb R$ for each $r\in\Bbb R$.
It follows that $\log\Delta_{\Cal R,\mu}$ is ergodic, i.e. $\Cal R$ is of type $III_1$.
Thus we obtain a new short proof of Lemma~1.3$(ii)$.
Moreover,
 for each   $t,r\in\Bbb R$, we have that $(tr,(1-t)r)$ is an essential value of
 the cocycle $t\alpha\times (1-t)\log\Delta_{\Cal R,\mu}$.
This yields  that $r$ is
 is an essential value of  the cocycle $\log\Delta_{\Cal R,\mu}-t\Lambda=t\alpha+(1-t)\log\Delta_{\Cal R,\mu}$ of $\Cal R$.
We thus obtain that the cocycle $\log\Delta_{\Cal R,\mu}-t\Lambda$  is ergodic for each real $t \in\Bbb R$.
\endremark

\remark{Remark 1.7}
The group $\Sigma_0$ of finite permutations of $\Bbb N$ acts naturally on $X$.
Denote by $S$ the $\Sigma_0$-orbit equivalence relation.
Of course, $\Cal S$ is a subrelation of $\Cal R$.
It follows from the proof of Proposition~1.5 (see the definition of $\phi$)  and Remark~1.6 that the restriction of the cocycle $\log\Delta_{\Cal R,\mu}-t\Lambda$  to $\Cal S$ is also ergodic for each real $t\in\Bbb R$.
\endremark

\head 2. Nonsingular endomorphisms
\endhead

Let $T$ be a countable-to-one $\mu$-nonsingular endomorphism of a $\sigma$-finite standard measure space $(X,\goth B,\mu)$.
The {\it $\mu$-nonsingularity} means that $\mu\circ T^{-1}\sim\mu$.
If for each $A\in\goth B$ of positive measure there is $n>0$ such that $\mu(T^{-n}A\cap A)>0$ then $T$ is called {\it conservative}.
If $T$ is not conservative then $T$ is called {\it dissipative}.
If $\mu(T^{-1}A\triangle A)=0$  implies $\mu(A)=0$ or $\mu(X\setminus A)=0$ then $T$ is called {\it ergodic}.
If $\bigcap_{n>0}T^{-n}\goth B=\{\emptyset,X\}$ (modulo the subsets of zero $\mu$-measure)
then $T$ is called {\it exact}.

The {\it orbit equivalence relation} $\Cal R_T$ of $T$ is given by the formula
$(x,y)\in \Cal R_T$ if and only if there are $n,m\ge 0$ such that $T^nx=T^my$.
Then $\Cal R_T$ is $\mu$-nonsingular.
We also consider a subrelation $\Cal S_T$ of $\Cal R_T$:  $(x,y)\in \Cal S_T$ if and only if there is $n\ge 0$ such that $T^nx=T^ny$.
Of course, $\Cal S_T$ is also $\mu$-nonsingular.
We recall some standard facts.

\proclaim{Lemma 2.1}
\roster
\item"$(i)$"
$\Cal R_T$ is $\mu$-ergodic if and only if $T$ is ergodic \cite{Haw}.
\item"$(ii)$"
$\Cal S_T$ is $\mu$-ergodic if and only if $T$ is exact \cite{Haw}.
\item"$(iii)$"
If  \,$T$ is invertible then $\Cal R_T$ is $\mu$-conservative if and  only if $T$ is conservative\footnote{This claim is no longer true in the general case where $T$ is non-invertible.
}.
\endroster
\endproclaim

From now on we will assume that $T$ is {\it aperiodic}, i.e. $\mu(\{x\in X\mid T^nx=x\})=0$ for each $n>0$.
Then for each  Borel function $\phi:X\to G$, there is a unique cocycle $\alpha_\phi$ of $\Cal R_T$ with values in $G$
such that $\alpha_\phi(x,Tx)=\phi(x)$.
By $T_\phi$ we denote the corresponding skew product transformation of $(X\times G,\mu\times\lambda_G)$:
$$
T_\phi(x,g)=(Tx,\phi(x)g).
$$
Of course, $T_\phi$ is a nonsingular countable-to-one endomorphism of  $(X\times G,\mu\times\lambda_G)$.
It is straightforward to verify that $\Cal R_{T_\phi}=\Cal R_T(\alpha_\phi)$
and $\Cal S_{T_\phi}=\Cal S_T(\alpha_\phi\restriction\Cal S_T)$.

If $T$ is invertible we denote by $\omega_{T,\mu}$
the {\it Radon-Nikodym derivative} $d\mu\circ T/d\mu:X\to\Bbb R^*_+$ of $T$.
Then, of course,  we have that
 $\alpha_{\omega_{T,\mu}}=\Delta_{\Cal R_T,\mu}$.
 This definition extends naturally to the general (non-invertible) case as follows.
Suppose that  $\mu$ is $\sigma$-finite on the $\sigma$-algebra $T^{-1}\goth B$.
Then by the {\it Radon-Nikodym derivative} of $T$ we mean the function $\omega_{T,\mu}:=(d\mu/d\mu\circ T^{-1})\circ T$.
However, in the non-invertible case we no longer have  that
 $\alpha_{\omega_{T,\mu}}=\Delta_{\Cal R_T,\mu}$.
 The endomorphism $T_{\omega_{T,\mu}}$ is called the {\it Maharam extension} of $T$.
 Choose a measure $\kappa$ on $\Bbb R^*_+$ equivalent to  Lebesgue measure such that
 $\kappa(aB)=a^{-1}\kappa(B)$ for each Borel subset $B\subset \Bbb R^*_+$ and $a\in \Bbb R^*_+$.
We will always assume that space $X\times\Bbb R^*_+$ of the Maharam extension is endowed with the measure $\mu\times\kappa$.
Then it is easy to see that $T_{\omega_{T,\mu}}$ preserves this measure.

We can associate a linear operator  $U_T$ in $L^2(X,\mu)$ to $T$ in the following way:
$$
U_Tf(x):=f(Tx)\sqrt{\omega_{T,\mu}(x)},\qquad x\in X.
$$
It is easy to see that $U_T$ is an isometry.
Hence $U_T^*U_T=I$.
If $T$ is invertible then $U_T$ is unitary and $U_T^n=U_{T^n}$ for all $n\in\Bbb Z$.\footnote{In the non-invertible case, in general, $U_T^n\ne U_{T^n}$ for $n>1$.}
A useful spectral condition for conservativeness of invertible endomorphisms was found in \cite{Ko2}.

\proclaim{Lemma 2.2 \cite{Ko2, Lemma~2.2}}
If $T$ is invertible and $\sum_{n\ge 0}\langle U_T^n 1, 1\rangle<\infty$ then $T$ is dissipative.
\endproclaim

If for each positive Borel function $f:X\to (0,+\infty)$,
$$
\sum_{n\ge 0}f(T^nx)\alpha_{\omega_{T,\mu}}(x,T^nx)=\infty \quad \text{a.e.}
$$
then $T$ is called {\it $\mu$-recurrent \cite{Si}}.

\proclaim{Lemma 2.3}
\roster
\item"$(i)$"
$T$ is $\mu$-recurrent if and only if $T_{\omega_{T,\mu}}$ is conservative \cite{Si}.
\item"$(ii)$"
If $T$ preserves $\mu$ then $T$ is $\mu$-recurrent if and only $T$ is conservative.
\item"$(iii)$"
 If $T$ is $\mu$-recurrent then $T$ is conservative \cite{Si}.
\item"$(iv)$"
If $T$ is invertible and  conservative then $T$ is
  $\mu$-recurrent.
\item"$(v)$"
If $T$ is non-invertible and conservative then there is a measure $\nu$ equivalent to $\mu$ such that $T$ is not $\nu$-recurrent \cite{EiSi}.
\item"$(vi)$" If $\mu(X)<\infty$ then $T$ is $\mu$-recurrent if and only if
$
\sum_{n\ge 0}\alpha_{\omega_{T,\mu}}(x,T^nx)=\infty \quad \text{a.e.}
$ \cite{Si}.
\endroster
\endproclaim

Suppose that $\mu$ is $\sigma$-finite on $T^{-1}\goth B$.
Then there is a standard $\sigma$-finite measure space $(\widetilde X,\widetilde{\goth B},\widetilde\mu)$, an invertible $\mu$-nonsingular transformation $\widetilde T$ of $\widetilde X$ and a Borel map $\pi:\widetilde X\to X$ such that the following are satisfied:
\roster
\item"$\bullet$"
$\widetilde\mu\circ \pi^{-1}=\mu$,
\item"$\bullet$"
$\pi\widetilde T=T\pi$,
\item"$\bullet$"
$\omega_{\widetilde T,\widetilde\mu}$ is $\pi^{-1}(\goth B)$-measurable,
\item"$\bullet$"
$\bigvee_{n>0} \widetilde T^n \pi^{-1}(\goth B)=\widetilde{\goth B}$
(mod $\widetilde\mu$).
\endroster
The dynamical system $(\widetilde X,\widetilde\mu,\widetilde T)$ is called {\it the natural extentsion} of $(X,\mu,T)$.
The natural extension exists and it is
unique up to a natural isomorphism  (see \cite{Si}, \cite{SiTh}).
In the case when $T$ preserves a probability measure, the natural extension of $T$ coincides with the well known Rokhlin's natural extension of $T$.

\example{Example 2.4 \cite{DaHa}}
Let $T$ be  a one-sided shift on $X=\{0,1\}^\Bbb N$ endowed with an infinite  product measure $\mu=\bigotimes _{i>0}\mu_i$, where $\mu_i$ is a distribution on $\{0,1\}$ such that $0<\mu_i(0)<1$ for each $i$.
Since $\mu\circ T^{-1}=\bigotimes _{i>0}\mu_{i+1}$, it follows from the
 Kakutani theorem on equivalence of infinite product measures \cite{Ka}
that $T$ is an endomorphism of $(X,\mu)$, i.e. $T$ is $\mu$-nonsingular, if and only if
$$
\sum_{i=1}^\infty
\bigg(\bigg(\sqrt{\mu_i(0)}-\sqrt{\mu_{i+1}(0)}\bigg)^2 +
\bigg(\sqrt{\mu_i(1)}-\sqrt{\mu_{i+1}(1)}\bigg)^2\bigg)<\infty.
\tag2-1
$$
We let $\widetilde X:=\{0,1\}^\Bbb Z$ and $\widetilde\mu:=\bigotimes_{n\in\Bbb Z}\widetilde\mu_n$, where $\widetilde\mu_n=\mu_1$ if $n\le 0$ and
$\widetilde\mu_n=\mu_n$ if $n\ge 1$.
Let $\widetilde T$ denote the two-sided shift on $\widetilde X$ and let $\pi:\widetilde X\to X$
denote the restriction map, i.e. $(\pi(x))_n=x_n$ for $n\ge 1$.
Then
$$
\frac{d\widetilde\mu\circ \widetilde T}{d\widetilde\mu}(x)=\prod_{i\in\Bbb Z}\frac{\widetilde\mu_{i-1}(x_i)}
{\widetilde\mu_{i}(x_i)}=
\prod_{i\ge 2}\frac{\widetilde\mu_{i-1}(x_i)}
{\widetilde\mu_{i}(x_i)},\qquad x\in\widetilde X.
$$
On the other hand,
$
\frac{d\mu}{d\mu\circ T^{-1}}(x)=\prod_{i\ge 1}\frac{\mu_i(x_i)}{\mu_{i+1}(x_i)}$.
Thus $\frac{d\widetilde\mu\circ \widetilde T}{d\widetilde\mu}(x)=\frac{d\mu}{d\mu\circ T^{-1}}(T\pi(x))$ for a.a. $x\in X$.
Hence $(\widetilde X,\widetilde\mu,\widetilde T)$ is the natural extension of $(X,\mu,T)$.
\endexample

A $\sigma$-finite measure $\nu$ on $(X,\goth B)$ is called {\it $T$-cohomologous} to $\mu$ if $\nu$ is equivalent to $\mu$ and the Radon-Nikodym derivative $d\nu/d\mu$ is measurable with respect to $T^{-1}\goth B$.

\proclaim{Lemma 2.5}
\roster
\item"$(i)$"
$(\widetilde X,\widetilde\mu,\widetilde T)$ is conservative if and only if $T$ is $\mu$-recurrent
\cite{SiTh}.
\item"$(ii)$"
If $\nu$ is $T$-cohomologous to $\mu$ then the natural extensions $(\widetilde X,\widetilde\mu,\widetilde T)$ and $(\widetilde X,\widetilde\nu,\widetilde T)$ of  $(X,\mu,T)$ and   $(X,\nu,T)$ respectively are isomorphic \cite{SiTh}.
\item"$(iii)$"
If $T$ is $\mu$-recurrent then  $T$ is ergodic if and only if $\widetilde T$ is ergodic.
\item"$(iv)$"
The Maharam extension  $\widetilde T_{\omega_{\widetilde T,\widetilde\mu}}$ of $\widetilde T$ is canonically isomorphic  to the natural extension $\widetilde{T_{\omega_{T,\mu}}}$ of the Maharam extension
of $T$.
\endroster
\endproclaim

Let $R$ be an invertible nonsingular transformation on a $\sigma$-finite standard measure space $(X,\goth B,\mu)$.
Then $R$ is said to be a {\it $K$-automorphism} (\cite{Pa}, \cite{SiTh}) if there is a $\sigma$-finite algebra $\goth F\subset \goth B$ such that
\roster
\item"$\bullet$"
$R^{-1}\goth F\subset\goth F$,
\item"$\bullet$"
the Radon-Nikodym derivative $\omega_{ R,\mu}$ is $\goth F$-measurable,
\item"$\bullet$"
$\bigvee_{n>0} R^n \goth F=\goth B$
(mod $\mu$),
\item"$\bullet$"
$\bigcap_{n>0} R^{-n} \goth F=\{\emptyset, X\}$ (mod $\mu$).
\endroster
In other words,  $R$ is a $K$-automorphisms if and only if it is the natural extension of an exact factor (semi-invariant  sub-$\sigma$-algebra) of $R$.

\proclaim{Lemma 2.6  (\cite{Pa}, \cite{SiTh})}
Each $K$-automorphism $R$ is either totally dissipative\footnote{We recall that an invertible transformation $R$ of a standard measure space $(Y,\nu)$ is called {\it totally dissipative} if there is a subset $Y_0\subset Y$ such that $R^nY_0\cap Y_0=\emptyset$ for each $n\in\Bbb N$ and $\bigsqcup_{n\in\Bbb Z }R^nY_0=Y$.} or conservative.
In the latter case  $R$  weakly mixing.
\endproclaim

The following example generalizes \cite{Ha, Theorem~1}\footnote{Hamachi considers only the case  $p=0.5$  in \cite{Ha}.}.
Our proof is shorter and more elementary.

\example{Example 2.7}
Let $\widetilde X=\{0,1\}^\Bbb Z$ and let $\mu=\bigotimes_{n\in\Bbb Z}\mu_n$ where $\mu_n(0):=p$ if $n\le 0$ and $\mu_n(0):= q$ if $n>0$ for some positive reals $p,q\in(0,1)$.
Let $\widetilde T$ denotes the two-sided shift on $\widetilde X$.
Of course, if $p=q$ then $\widetilde T$ is conservative
because $\widetilde T$ is a probability preserving Bernoulli shift.
We now show that if $p\ne q$ then $\widetilde T$ is dissipative.
Without loss of generality we may assume that $p<q$.
For $\mu$-a.a. $x\in \widetilde X$, we have
$$
\frac{d\mu\circ \widetilde T}{d\mu}(x)=\prod_{n\in\Bbb Z}\frac{\mu_{n-1}(x_{n})}{\mu_{n}(x_{n})}=
\frac{\mu_0(x_1)}{\mu_1(x_1)}.
$$
Therefore for each $n>0$,
$$
\frac{d\mu\circ \widetilde T^n}{d\mu}(x)=\frac{\mu_0(x_{1})\cdots\mu_0(x_n)}{\mu_1(x_1)\cdots\mu_1(x_n)}=\left(\frac pq\right)^n\left(\frac{q(1-p)}{p(1-q)}\right)^{x_1+\cdots+x_n}
\tag2-2
$$
at a.e. $x$.
Since $p^q(1-p)^{1-q}<q^q(1-q)^{1-q}$, there is $\epsilon>0$
such that
$$
\delta:=\left(\frac pq\right)^{q-\epsilon} \left(\frac {1-p}{1-q}\right)^{1-q+\epsilon}<1.
\tag2-3
$$
It follows from the  individual ergodic theorem (for the one-sided shift) that if $n$ is sufficiently large then
$\frac{x_1+\cdots+x_n}n\le 1-q+\epsilon$ for a.e. $x$.
It follows from \thetag{2-2} and \thetag{2-3} that
$
\frac{d\mu\circ\widetilde T^n}{d\mu}(x)\le\delta^n.
$
Therefore the series $\sum_{n\ge 1}\frac{d\mu\circ \widetilde T^n}{d\mu}(x)$ converges at a.e. $x$.
Hence $\widetilde T$ is not $\mu$-recurrent.
It follows from  Lemma~2.3$(iv)$ that $\widetilde T$ is dissipative.
Moreover, in view of Example~2.4, $\widetilde T$ is a natural extension of an exact endomorphism.
Hence, $\widetilde T$ is a $K$-automorphism.
 By  Lemma~2.6, $\widetilde T$
is totally dissipative.
\endexample

\head 3. Nonsingular Bernoulli shifts
\endhead

Throughout this section we will use the notation introduced in Example~2.4.
Thus  $(\widetilde X,\widetilde \mu,\widetilde T)$  stands for the natural extension of the one-sided $\mu$-nonsingular Bernoulli shift  $(X,\mu,T)$.
In particular, \thetag{2-1} holds.
It is easy to verify that
 $\Cal S_T$ is the tail equivalence relation on $X$.
 Hence $\Cal S_T$ is ergodic by Kolmogorov's zero-one law.
By Lemma~2.1$(ii)$, $T$ is exact.
Therefore $\widetilde T$ is a $K$-automorphism.

\proclaim{Theorem 3.1}
\roster
\item"$(i)$"
If $\sum_{n=1}^\infty(\mu_n(0)-\mu_1(0))^2<\infty$ then 
$\widetilde\mu$ is equivalent to $\bigotimes_{n=-\infty}^\infty\mu_1$
\footnote{Thus, in this case, $\widetilde T$ is isomorphic to the measure preserving Bernoulli shift on  $(\widetilde X,\bigotimes_{n=-\infty}^\infty\mu_1)$.}.
\item"$(ii)$"
If $\sum_{n=1}^\infty(\mu_n(0)-\mu_1(0))^2=\infty$ and $\widetilde T$ is conservative  then the Maharam extension  $\widetilde T_{\omega_{\widetilde T,\widetilde\mu}}$ of $\widetilde T$ is a weakly mixing  $K$-automorphism.
\endroster
\endproclaim

\demo{Proof}
$(i)$
It follows from Lemma~1.3$(i)$ that $\mu$ is equivalent to the infinite product $\nu:=\mu_1\times\mu_1\times\cdots$.
 Moreover,  it is easy to see that the Radon-Nikodym derivative $d\mu/d\nu$ does not depend on $x_1$.
 Hence by Lemma~2.5$(ii)$, the natural extension of $(T,\nu)$ is isomorphic to $\widetilde T$.
 On the other hand, according to Example~2.4, the natural extension of $(T,\nu)$  is the  2-sided Bernoulli shift on
 $(\widetilde X, \bigotimes_{n\in\Bbb Z}\mu_1)$.

$(ii)$
By the Maharam theorem, since $\widetilde T$ is conservative,
$\widetilde T_{\omega_{\widetilde T,\widetilde\mu}}$ is also conservative.
Hence in view of Lemma~2.6, if $\widetilde T_{\omega_{\widetilde T,\widetilde\mu}}$ is a $K$-automorphism then it is weakly mixing.
Thus it suffices to prove that $\widetilde T_{\omega_{\widetilde T,\widetilde\mu}}$ is a $K$-automorphism.
We will proceed in several steps.
Let $\Cal L$ stand for the set of limit points of the sequence $(\mu_n(0))_{n\ge 1}$.

{\it Claim A.}   $\Cal L\ni\mu_1(0)$.

Indeed, if  $\Cal L\not\ni\mu_1(0)$ then there are $\delta>0$ and $N>0$ such that $|\mu_1(0)-\mu_n(0)|>\delta$ for each  $n\ge N$.
Therefore we have that
$$
\align
\int_{\widetilde X}\sqrt{\frac{d\widetilde\mu\circ\widetilde T^n}{d\widetilde \mu}(x)}\,d\widetilde\mu(x)&=
\prod_{k\in\Bbb Z}\int \sqrt{\frac{d\mu_{k-n}}{d \mu_k}(x_k)}\,d\mu_k(x_k)\\
&=
\prod_{k>0}(\sqrt{\mu_{k-n}(0)\mu_k(0)}+\sqrt{\mu_{k-n}(1)\mu_k(1)}\,)\\
&\le
\prod_{k=N}^n(\sqrt{\mu_{k-n}(0)\mu_k(0)}+\sqrt{\mu_{k-n}(1)\mu_k(1)}\,)\\
&=
\prod_{k=N}^n(\sqrt{\mu_{1}(0)\mu_k(0)}+\sqrt{\mu_{1}(1)\mu_k(1)}\,).
\endalign
$$
Since $\eta:=\sup\{\sqrt{\mu_1(0)t}+\sqrt{\mu_1(1)(1-t)}\mid t\in(0,1), |t-\mu_1(0)|>\delta\}<1$, we obtain that $\langle U_{\widetilde T}^{n}1,1\rangle\le\eta^{n-N}$.
It now follows from Lemma~2.2 that $\widetilde T$ is dissipative.
This contradicts to the condition of the theorem.
Hence  $\mu_1(0)\in \Cal L$, as claimed.

 {\it Claim B.}
  If $\Cal L=\{\mu_1(0)\}$
 then $\widetilde T_{\omega_{\widetilde T,\widetilde\mu}}$ is a $K$-automorphism.

 By Lemma~2.5$(iv)$, $\widetilde T_{\omega_{\widetilde T,\widetilde\mu}}$ is isomorphic to the natural extension
 $\widetilde{T_{\omega_{T,\mu}}}$ of the endomorphism
$T_{\omega_{T,\mu}}$.
Hence  $\widetilde T_{\omega_{\widetilde T,\widetilde\mu}}$ is a $K$-automorphism if    $T_{\omega_{T,\mu}}$ is exact.
By Lemma~2.1$(ii)$, $T_{\omega_{T,\mu}}$ is exact if and only if the equivalence relation $\Cal S_{T_{\omega_{T,\mu}}}$ is ergodic.
We note that $\Cal S_{T_{\omega_{T,\mu}}}=\Cal S_T(\alpha_{\omega_{T,\mu}}\restriction\Cal S_T)$.
Thus it suffices  to show that   the cocycle  $\alpha_{\omega_{T,\mu}}$ restricted to $\Cal S_T$ is ergodic.
Take $(x,y)\in S_T$.
  Then there is $n>0$ such that $T^nx=T^ny$, i.e. $x_{i}=y_i$ if $i>n$.
  Let $\widetilde x,\widetilde y\in\widetilde X$ be such that $\pi(\widetilde x)=x$ and
 $\pi(\widetilde y)=y$.
Since
$$
\alpha_{\omega_{T,\mu}}(x,T^nx)=\omega_{T,\mu}(x)\cdots\omega_{T,\mu}(T^{n-1}x)=
\omega_{\widetilde T,\widetilde\mu}(\widetilde x)\cdots\omega_{\widetilde T,\widetilde\mu}(\widetilde T^{n-1}\widetilde x)=\frac{d\widetilde\mu\circ\widetilde T^n}{d\widetilde\mu}(\widetilde x),
$$
it follows that
 $$
 \alpha_{\omega_{T,\mu}}(x,y)
=\frac{ \alpha_{\omega_{T,\mu}}(x,T^nx)}{ \alpha_{\omega_{T,\mu}}(y,T^ny)}
=
 \frac{\frac{d \widetilde\mu\circ \widetilde T^n}{d\widetilde\mu}(\widetilde x)}{\frac{d \widetilde\mu\circ \widetilde T^n}{d\widetilde\mu}(\widetilde y)}.
$$
This yields
$$
\align
 \alpha_{\omega_{T,\mu}}(x,y)
&=
 \frac{\prod_{i=1}^n\frac{\mu_1(x_i)}{\mu_i(x_i)}\prod_{i=n+1}^\infty \frac{\mu_{i-n}(x_i)}{\mu_i(x_i)}}
 {\prod_{i=1}^n\frac{\mu_1(y_i)}{\mu_i(y_i)}\prod_{i=n+1}^\infty \frac{\mu_{i-n}(y_i)}{\mu_i(y_i)}}\\
&
=
\prod_{i=1}^n\frac{\mu_i(y_i)}{\mu_i(x_i)}\prod_{i=1}^n\frac{\mu_1(x_i)}{\mu_1(y_i)}.
\endalign
 $$
Thus we obtain that
$
 \alpha_{\omega_{T,\mu}}(x,y)=\Delta_{\Cal S_T,\mu}(x,y)/\delta(x,y),
$
where $\delta$ is a cocycle of $\Cal S_T$ is given by
$\delta(x,y)=\prod_{i=1}^\infty\mu_1(y_i)/\mu_1(x_i)$.
It follows from~Proposition~1.4 that
 $\alpha_{\omega_{T,\mu}}$ restricted to $\Cal S_T$ is ergodic.

{\it Claim C.} If $\mu_1(0)\in\Cal L$ but $\Cal L\ne\{\mu_1(0)\}$ then
$\widetilde T_{\omega_{\widetilde T,\widetilde\mu}}$ is a $K$-automorphism.

As we have proved in Claim~B, it suffices to show that
 $\alpha_{\omega_{T,\mu}}\restriction\Cal S_T$ is ergodic.
 Since   \thetag{2-1} is satisfied,
there is a segment $[\alpha,\beta]\subset\Cal L$ such that $\mu_1(0)\in [\alpha,\beta]$.
 Then it is easy to see that that we can find an infinite subset $I\subset\Bbb N$,
a real  $\lambda\in[\alpha,\beta]$ and a sequence $\epsilon_i\to 0$ such that the complement $\Bbb N\setminus I$ is infinite,
 $\mu_i(0)=\lambda e^{\epsilon_i}$ if $i\in I$ and $\sum_{i\in I}\epsilon_i^2=\infty$.
 Since $\Cal S_T$ is naturally isomorphic to the direct product of the tail equivalence relations $\Cal S^1$ and $\Cal S^2$ on $(\{0,1\}^I,\bigotimes_{i\in I}\mu_i)$ and $(\{0,1\}^{\Bbb N\setminus I},\bigotimes_{i\in \Bbb N\setminus I}\mu_i)$ respectively and  the restriction of $ \alpha_{\omega_{T,\mu}}$ to $\Cal S^1$ is ergodic by Claim~B, it follows that $ \alpha_{\omega_{T,\mu}}\restriction\Cal S_T$ is also ergodic.
\qed
\enddemo

\remark{Remark 3.2}
\roster
\item"(i)"
We note that if  the Bernoulli shift $(X,\mu, T)$ is {\it equilibrial}, i.e. $\mu_1(0)=\mu_1(1)$ then  (and only in this case)
$ \alpha_{\omega_{T,\mu}}\restriction \Cal S_T=\Delta_{\Cal S_T,\mu}$.
Therefore in this case to prove Theorem~3.1 it would suffice to apply the well known Lemma~1.3$(ii)$ instead of Propostion~1.4.
That was done in \cite{Ko2}.
\item"(ii)"
The Krengel entropy of $\widetilde T$ is infinite \cite{SiTh}.
\endroster
\endremark

The next statement follows  immediately from  Theorem~3.1.
\proclaim{Corollary 3.3}
If  $\sum_{n=1}^\infty(\mu_1(0)-\mu_n(0))^2<\infty$
then $\widetilde T$ is of type $II_1$.
If \,$\sum_{n=1}^\infty(\mu_1(0)-\mu_n(0))^2=\infty$
and  $\widetilde T$ is conservative then
  $\widetilde T$ is
of type $III_1$.
\endproclaim

\head 4. Forcing conservativeness of Bernoulli shifts
\endhead

Let $X=\{0,1\}^\Bbb Z$ and
  let
$$
\Cal A_n:=\{[a_n,\dots ,a_n]_{-n}^n\mid a_{-n},\dots, a_n\in\{0,1\}\}.
 $$
Denote by $T$ the 2-sided shift on $X$.
We now state without proof  a standard approximation result (cf. Lemma~1.1).

\proclaim{Lemma 4.1}
Let $\mu$ be a probability measure on $X$ and let $T$ be $\mu$-nonsingular.
If for each $n>0$ and $A\in\Cal A_n$, there are a subset $A_0\subset A$ and a one-to-one transformation $\tau_A:A_0\to A$ such that $\mu(A_0)>0.9\mu(A)$ and $\tau_A x\in\{T^nx\mid n>0\}$ and
$$
\left|\log\left(\frac{d\mu\circ\tau_A}{d\mu}(x)\right)\right|<0.001
$$
for a.e. $x\in A_0$ then $T$ is $\mu$-conservative.
\endproclaim

 Fix $\lambda\in (0,1]$ and  a sequence $(\epsilon_n)_{n\in\Bbb Z}$ of reals such that
$\epsilon_n=0$ if   $n\le 1$ and
$$
\lim_{n\to\infty}\epsilon_n=0.\tag4-1
$$
We define a measure $\mu$ on $X$ by setting
 $$
\mu=\bigotimes_{n\in\Bbb Z}\mu_n, \text{ where }
\mu_n(0)=\frac 1{1+\lambda e^{\epsilon_n}},
\mu_n(1)=\frac {\lambda e^{\epsilon_n}}{1+\lambda e^{\epsilon_n}}\text{ for  each $n\in\Bbb Z$.}\tag4-2
$$
It follows from  Kakutani's theorem (see \thetag{2-1} and \cite{Ka}), that $T$  is $\mu$-nonsingular if and only if
$$
\sum_{n=1}^\infty|\epsilon_{n+1}-\epsilon_n|^2<\infty.\tag4-3
$$
According to Theorem~3.1, if $T$ is conservative and
$$
\sum_{n=1}^\infty\epsilon_n^2=\infty\tag4-4
$$
then $T$ is ergodic of type $III_1$.
 Thus our purpose is to construct $(\epsilon_n)_{n\ge 1}$ such that ~\thetag{4-1}, \thetag{4-3} and \thetag{4-4} are satisfied and $T$ is conservative.
 We will do this inductively.
 Each step of the inductive construction will consist of two semi-steps.
On the first semi-step we ``do conservativeness'' of  $T$ partly, on $\Cal A_n$.
 On the second semi-step we ``satisfy partly'' \thetag{4-1}, \thetag{4-3} and \thetag{4-4}.
The only additional problem
 is that the second semi-step of the $m$-th step will affect the the property of partial conservativeness  achieved on the $n$-th steps for $n<m$.
 Thus we have to control that the total contribution of the subsequent steps $(m>n)$ into the partial conservativeness of $T$ on $\Cal A_n$ is ``small''.

 Fix a sequence $(\eta_n)_{n=1}^\infty$ of positive reals such that $\eta_n\to 0$ as $n\to\infty$.
Suppose that we have already defined $\epsilon_1,\dots,\epsilon_{L_{n-1}}$.
We now let
$$
\mu^{(n)}:=\bigotimes_{i\le L_{n-1}}\mu_i\otimes\bigotimes_{i>L_{n-1}}\mu_1
$$
Since $\mu^{(n)}$ is equivalent to the infinite product $\bigotimes_{i\in\Bbb Z}\mu_1$,
it follows that
$T$ is $\mu^{(n)}$-nonsingular
and $\mu^{(n)}$-conservative.
Hence applying the standard exhaustion argument we can find
 for  each cylinder $A\in\Cal A_{L_{n-1}}$, positive integers $p_1,\dots, p_{m}$
and  pairwise disjoint cylinders $B_1,\dots, B_m\subset A$
such that
$$
T^{p_i}B_i\subset A,\ T^{p_i}B_i\cap T^{p_j}B_j=\emptyset\ \text{if }1\le i\ne j\le m
\text{ and}\
\mu^{(n)}\left(\bigsqcup_{i=1}^mB_i\right)>0.9\mu^{(n)}(A).
$$
We now define a map $\tau_A:\bigsqcup_{i=1}^mB_i\to A$ by setting
$
\tau_A x:=T^{p_i}x
$
if $x\in B_i$.
Then $\tau_A$ is one-to-one and $ \frac{d\mu^{(n)}\circ \tau_A}{d\mu^{(n)}}(x)=1$ for all $x\in \bigsqcup_{i=1}^mB_i $.
Choose $\ell_n>L_{n-1}$ large so that  $B_i$ and $T^{p_i}B_i$ are unions of cylinders from  $\Cal A_{\ell_n}$ for each $i$.
It follows that $p_i<\ell_n$ for $i=1,\dots,m$.
We now set $\epsilon_j:=0$ if $L_{n-1}<j\le \ell_n$.

Now we choose an integer $L_n>2\ell_n$ and reals $\epsilon_{\ell_n+1}\ge\cdots\ge\epsilon_{L_n}$ so
that
$$
\sum_{j=\ell_n+1}^{2\ell_n}\epsilon_j<\eta_n,\  \sum_{j=\ell_n+1}^{L_n}(\epsilon_{j-1}-\epsilon_j)^2<\eta_n\text{ \ and } \sum_{j=\ell_n+1}^{L_n}\epsilon_j^2>1.\tag4-5
$$
Thus we defined $\epsilon_1,\dots,\epsilon_{L_n}$.
Continuing this construction process infinitely many times, we obtain an increasing sequence $\ell_1<L_1<\ell_2<L_2<\cdots$
and  a sequence $(\epsilon_n)_{n\ge 1}$.
Moreover,
for each cylinder $A\in\bigcup_{n\ge 1}\Cal A_{L_n}$, we have a map $\tau_A$ satisfying the properties listed above.

\proclaim{Theorem 4.2}
Let $\lambda\in(0,1)$ and let $(\epsilon_n)_{n>0}$ be a sequence of nonnegative reals defined via the aforementioned inductive procedure.
Define $\mu$  by \thetag{4-2} utilizing $\lambda$ and $(\epsilon_n)_{n>0}$.
Then the  2-sided shift $T$ on  $(X,\mu)$ is $\mu$-nonsingular, $\mu$-conservative and  of type $III_1$.
\endproclaim

\demo{Proof}
Since \thetag{4-5} implies \thetag{4-1}, \thetag{4-3} and \thetag{4-4},
$T$ is $\mu$-nonsingular and 
of type $III_1$ whenever it is $\mu$-conservative (in view of Corollary~3.3).
Thus it suffices  to verify that $T$ is $\mu$-conservative.
For that we will apply~Lemma~4.1.

We first note that $\mu(C)=\mu^{(n)}(C)$ for each cylinder $C\in\Cal A_{\ell_n}$, $n>0$.
Take a cylinder $A\in\Cal A_{L_{n-1}}$.
By the definition of $\mu$, there are mutually disjoint cylinders $B_1,\dots, B_m\in\Cal A_{\ell_n}$, positive integers  $p_1,\dots,p_m<\ell_n$ and a one-to-one map
$\tau_A:\bigsqcup_{i=1}^m B_i\to A$
such that
$
\mu^{(n)}\left(\bigsqcup_{i=1}^m B_i\right)>0.9\mu^{(n)}(A)$,
 $\tau_Ax=T^{p_i}x$  and $\frac{d\mu^{(n)}\circ T^{p_i}}{d\mu^{(n)}}(x)=1$ for each $x\in B_i$.
It follows that
$\mu\left(\bigsqcup_{i=1}^m B_i\right)>0.9\mu(A)$ and
$$
\frac{d\mu\circ T^{p_i}}{d\mu}(x)=\frac{d\mu}{d\mu^{(n)}}(T^{p_i}x)\frac{d\mu^{(n)}}
{d\mu}(x)
\quad
 \text{for each $x\in B_i$.}
$$
Since
$
\frac{d\mu}{d\mu^{(n)}}(x)=\prod_{k>n}\prod_{j=\ell_k+1}^{L_k}
\frac{\mu_j(x_j)} {\mu_1(x_j)}
$
for each $x\in X$, we obtain  that
$$
\frac{d\mu\circ T^{p_i}}{d\mu}(x)=\prod_{k>n}\prod_{j=\ell_k+1}^{L_k}
\frac{\mu_j(x_{j+p_i})\mu_1(x_j)} {\mu_j(x_j)\mu_1(x_{j+p_i})}
=\prod_{k>n}\prod_{j=\ell_k+1}^{L_k+p_i}
\frac{\mu_{j-p_i}(x_{j})} {\mu_j(x_j)}\tag4-6
$$
for each $x\in B_i$.
Given $0<a<b$ and $p>0$, we have
$$
\prod_{j= a}^{b}\frac{\mu_{j-p}(x_j)}{\mu_j(x_j)}=
\prod_{j=a}^b\frac{(\lambda e^{\epsilon_{j-p}})^{x_j}}{1+\lambda e^{\epsilon_{j-p}}}\cdot
\frac{1+\lambda e^{\epsilon_j}}{(\lambda e^{\epsilon_j})^{x_j}}
=
\prod_{j=a}^b
e^{x_j(\epsilon_{j-p}-\epsilon_j)}
\frac{1+\lambda e^{\epsilon_j}}{1+\lambda e^{\epsilon_{j-p}}}.
\tag4-7
$$
 We  note that
 $\epsilon_{j-p}\ge\epsilon_j$ if $L_n+p\ge j\ge \ell_n+1+p$ and $p<\ell_n$.
 Hence
  $$
  \align
  - \sum_{j=\ell_n+1}^{\ell_n+p}\epsilon_j \le &
\sum_{j=\ell_n+1}^{L_n+p}x_j(\epsilon_{j-p}-\epsilon_j)\le\sum_{j=\ell_n+1+p}^{L_n+p}(\epsilon_{j-p}-\epsilon_j)=\sum_{j=\ell_n+1}^{\ell_n+p}\epsilon_j.
\tag4-8
  \endalign
  $$
On the other hand, we see that
$$
\prod_{j=\ell_n+1}^{L_n+p}
\frac{1+\lambda e^{\epsilon_j}}{1+\lambda e^{\epsilon_{j-p}}}=
\frac{\prod_{j=L_n+1} ^{L_n+p}(1+\lambda e^{\epsilon_j})}{\prod_{j=\ell_n+1-p}^{\ell_n}(1+\lambda e^{\epsilon_j})}=\frac{(1+\lambda)^p}{(1+\lambda)^p}=1.
$$
This, \thetag{4-6}--\thetag{4-8} and \thetag{4-5} yield that for each $x\in\bigsqcup_{i=1}^mB_i$,
$$
\log\left|\frac{d\mu\circ \tau_A}{d\mu}(x)\right|=\log\left|\frac{d\mu\circ T^{p_i}}{d\mu}(x)\right|\le
\sum_{j=\ell_n+1}^{\ell_n+p_i}\epsilon_j
\le\eta_n
$$
where $i$ is chosen so that $x\in B_i$.
Hence $T$ is $\mu$-conservative by Lemma~4.1.
\qed
\enddemo

\head 5.  Markov measures
\endhead

Let $X=\{0,1\}^{\Bbb N}$.
Given a distribution $\lambda$ on $\{0,1\}$ and a sequence $\boldsymbol P:=(P^{(n)})_{n=1}^\infty$ of stochastic $2\times 2$ matrices $P^{(n)}=(P^{(n)}_{i,j})_{i,j=0,1}$, we define a Borel measure $\mu$ on $X$  by setting
$$
\mu([a_1,\dots,a_k]_1^k):=\lambda(a_1)P^{(1)}_{a_1,a_2}P^{(2)}_{a_2,a_3}\cdots P^{(k-1)}_{a_{k-1},a_k}.
$$
It is called the {\it Markov measure} determined by the pair $(\lambda,\boldsymbol P)$.
We say that $\mu$ is {\it non-degenerate} if  $\lambda(a)>0$ and $P^{(n)}_{a,b}>0$ for all $a,b\in\{0,1\}$ and $n>0$.
It is easy to see that $\mu$ is {\it non-atomic} if and only if  $\prod_{n=1}^\infty P^{(n)}_{x_n,x_{n+1}}=0$ for each $x\in X$.
Let $\Cal R$ denote the tail equivalence relation on $X$.
If $\mu$ is non-degenerate then  $\Cal R$ is $\mu$-nonsingular.
Indeed, it is straightforward  to verify that
$$
\Delta_{\Cal R,\mu}(x,y)=\frac{\lambda(y_1)}{\lambda(x_1)}\prod_{j=1}^\infty\frac{P^{(j)}_{y_j,y_{j+1}}}{P^{(j)}_{x_j,x_{j+1}}},
$$
where the product is, in fact, finite because $x$ and $y$ are $\Cal R$-equivalent.
The following theorem is a generalization of the well known Kakutani theorem on equivalence of infinite product measures \cite{Ka}.

\proclaim{Lemma 5.1 \cite{Lo}}
Let $\mu$ and $\nu$ be two Markov measures on $X$ determined by pairs $(\lambda,\boldsymbol P)$ and $(\kappa,\boldsymbol Q)$ respectively.
Let the $\sigma$-algebra of $\Cal R$-invariant subsets be trivial $\pmod \nu$ and let $\mu$ be non-degenerate.
Then 
 $\nu$ is absolutely continuous with respect to $\mu$  if and only if
$$
\int_X\sqrt{\frac{\kappa(x_1)Q^{(1)}_{x_1,x_2}Q^{(2)}_{x_2,x_3}\cdots Q^{(k-1)}_{x_{k-1},x_k}}{\lambda(x_1)P^{(1)}_{x_1,x_2}P^{(2)}_{x_2,x_3}\cdots P^{(k-1)}_{x_{k-1},x_k}}}\,d\mu(x)\not\to 0
\tag5-1
$$
as $k\to\infty$.
\endproclaim

A natural question arises: under which conditions on $\mu$ the $\sigma$-algebra of Borel $\Cal R$-invariant subsets is trivial $\pmod \mu$?
To answer it, we  state a theorem from \cite{BrDo} which is an analog of Kolmogorov zero-one law for the general probability measures on $X$.
For that we need a piece of notation.
Denote by $\goth B_n$ the (finite) $\sigma$-algebra generated by cylinders $[a_1,\dots,a_n]_{1}^n$, $a_1,\dots,a_n\in\{0,1\}$.
Denote by $\goth B^n$ the smallest $\sigma$-algebra  such that the maps $X\ni x\mapsto x_k\in\{0,1\}$, $k>n$, are all measurable.
Let $\mu$ be a probability measure on $X$.
For each $n>0$, we denote by $\widehat\mu_n$ the following probability measure on $X$ determined by
$$
\widehat\mu_n(A\cap B)=\mu(A)\mu(B),\qquad A\in\goth B_n, B\in\goth B^n.
$$
This measure is equivalent to $\mu$.
We let  $r_n:=d\widehat\mu_n/d\mu$.

\proclaim{Lemma 5.2 \cite{BrDo}}
Let $\mu$ be a probability measure on $X$.
The $\sigma$-algebra of Borel $\Cal R$-invariant subsets is trivial (mod $\mu$) if and only
if $E(r_n\mid\goth B_n\vee \goth B^l)\to 1$ in measure $\mu$ as $l\to\infty$ for each $n\ge 1$.
\endproclaim

Let $\mu$ be a  Markov measure determined by $(\lambda,\boldsymbol P)$.
Given $m>n>0$, we denote by $P^{(n,m)}$ the matrix product  $P^{(n)}P^{(n+1)}\cdots P^{(m)}$.
Of course, $P^{(n,m)}$ is also a stochastic matrix.

\proclaim{Lemma 5.3} If $\mu$ is non-degenerate and non-atomic then for each $n>0$ and $a,b,c\in\{0,1\}$, there exists a limit
$$
\lim_{m\to\infty}\frac{ P^{(n,m)}_{a,c}}{P^{(n,m)}_{b,c}}= 1.
$$
\endproclaim

\demo{Proof} 
Of course, it suffices to verify only the case $a=c=0$ and  $b=1$.
To this end, we  first note that for each stochastic $2\times 2$ matrix $A=\pmatrix a_{0,0}&a_{0,1}\\a_{1,0}&a_{1,1}\endpmatrix$ with non-zero entries, $|\det A|=|a_{0,0}-a_{1,0}|<1$.
Hence either $\prod_{m=n}^\infty|\det P^{(m)}|=0$ for  each $n>0$ or   $\prod_{m=n}^\infty|\det P^{(m)}|>0$ for some $n>0$.
In the latter case, we obtain that  
$$
\prod_{m=n}^\infty\max(P^{(m)}_{0,0},P^{(m)}_{1,0})\ge\prod_{m=n}^\infty|P^{(m)}_{0,0}-P^{(m)}_{1,0}|=\prod_{m=n}^\infty|\det P^{(m)}|>0.
$$
This contradicts to the assumption that $\mu$ is non-atomic.
Hence  
$$
\lim_{m\to\infty}|P^{(n,m)}_{0,0}-P^{(n,m)}_{1,0}|=\lim_{m\to\infty} |\det P^{(n,m)}|=\prod_{m=n}^\infty|\det P^{(m)}|=0\tag5-2
$$
 for  each $n>0$.
Moreover, the sequence $(\min(P^{(n,m)}_{0,0},P^{(n,m)}_{1,0}))_{m>n}$ increases and the sequence $(\max(P^{n,m}_{0,0},P^{n,m}_{1,0}))_{m>n}$ decreases as $m\to\infty$ for each $n$.
Hence \thetag{5-2} implies that $\frac{P^{(n,m)}_{0,0}}{P^{(n,m)}_{1,0}}\to 1$ as $m\to\infty$.
\qed
\enddemo

\proclaim{Theorem 5.4} Let $\mu$ be a non-degenerate and non-atomic Markov measure determined by $(\lambda,\boldsymbol P)$.
 Then $\Cal R$ is $\mu$-ergodic.
\endproclaim
\demo{Proof}
Let $\mu=\int\mu_x\,d\nu(x)$ stand for the disintegration of  $\mu$ with respect to the restriction $\nu$ of $\mu$ to $\goth B_n\vee \goth B^{l-1}$.
For each finite sequence $a_1,\dots,a_m\in\{0,1\}$ with $m>l$, we have
$\mu_x([a_1,\dots,a_m]_1^m)=0$ if $a_j\ne x_j$
 for some $j\in\{1,\dots,n\}\cup\{l,\dots m\}$ and 
$$
\align
\mu_x([a_1,\dots,a_m]_1^m)&=\frac{\mu([x_1,\dots,x_n,a_{n+1},\dots a_{l-1},x_l,\dots,x_m]_l^m)}{\mu([x_1,\dots,x_n]_1^n\cap[x_l,\dots,x_m]_l^m)}\\
&=\frac{P^{(n)}_{x_n,a_{n+1}}P^{(n+1)}_{a_{n+1},a_{n+2}}\cdots
P^{(l-1)}_{a_{l-1},x_l}}
{\sum_{b_{n+1},\dots,b_{l-1}=0}^1{P^{(n)}_{x_n,b_{n+1}}P^{(n+1)}_{b_{n+1},b_{n+2}}
\cdots
P^{(l-1)}_{b_{l-1},x_l}}}\cdot\frac{\mu([x_n]_n)}{\mu([x_n]_n)}\\
&=\frac{\mu([x_n,a_{n+1},\dots,a_{l-1},x_{l}]_n^{l})}{\mu([x_n]_n\cap[x_{l}]_{l})}
\endalign
$$
otherwise.
We also have that
$$
\align
\frac{d\widehat\mu_n}{d\mu}(x)&
=
\lim_{m\to\infty}\frac{\widehat\mu_n([x_1,\dots,x_m]_1^m)}{\mu([x_1,\dots,x_m]_1^m)}\\
&=\lim_{m\to\infty}\frac{\mu([x_1,\dots,x_n]_1^n)\mu([x_{n+1},\dots,x_m]_{n+1}^m)}{\mu([x_1,\dots,x_m]_1^m)}
\\
&=\lim_{m\to\infty}\frac{\mu([x_n]_n)\mu([x_{n+1}]_{n+1})}{\mu([x_n,x_{n+1}]_n^{n+1})}\\
&=\frac{\mu([x_{n+1}]_{n+1})}{P^{(n)}_{x_n,x_{n+1}}}.
\endalign
$$
This yields
$$
\align
E(r_n\mid\goth B_n\vee \goth B^{l-1})(x)
&=\int\frac{\mu([t_{n+1}]_{n+1})}{P^{(n)}_{t_n,t_{n+1}}}\,d\mu_x(t)\\
&=\sum_{i=0}^1\mu([i]_{n+1})\frac{\mu([x_n,i]_n^{n+1}\cap[x_{l}]_{l})}
{P^{(n)}_{x_n,i}\mu([x_n]_n\cap[x_{l}]_{l})}\\
&=\sum_{i=0}^1\mu([i]_{n+1})\frac{P_{i,x_l}^{(n+1,l-1)}}{\sum_{j=0}^1P^{(n)}_{x_n,j}P^{(n+1,l-1)}_{j,x_l}}.
\endalign
$$
It follows from Lemma~5.3 and  the condition of the theorem that $E(r_n\mid\goth B_n\vee \goth B^{l-1})(x)\to 1$ at a.e. $x\in X$.
By Lemma~5.2,  $\Cal R$ is $\mu$-ergodic.
\qed
\enddemo

\proclaim{Corollary 5.5} Let $\mu$ and $\nu$ be two non-degenerate and non-atomic Markov measures on $X$ determined by pairs $(\lambda,\boldsymbol P)$ and $(\kappa,\boldsymbol Q)$ respectively.
Then 
 $\nu$ is equivalent to $\mu$  if and only if
 $$
 \lim_{n\to\infty}
 \sum_{a_1=0}^1\sqrt{\kappa(a_1)\lambda(a_1)}\prod_{i=1}^{n-1}\sum_{a_{i+1}=0}^1\sqrt{Q^{(i)}_{a_{i},a_{i+1}}P^{(i)}_{a_{i},a_{i+1}}}\ne 0.
 \tag5-3
 $$
 Moreover, in this case,
$$
\frac{d\nu}{d\mu}(x)=\frac{\kappa(x_1)}{\lambda(x_1)}\prod_{j=1}^\infty\frac{Q^{(j)}_{x_j,x_{j+1}}}{P^{(j)}_{x_j,x_{j+1}}}.
$$
If $\inf_{n>0}\min_{i,j}P^{(n)}_{i,j}>0$ and $\inf_{n>0}\min_{i,j}Q^{(n)}_{i,j}>0$  then $\mu$ and $\nu$ are equivalent if and only if $\sum_{n=1}^\infty(P^{(n)}_{i,j}-Q^{(n)}_{i,j})^2<\infty$ for each $i,j\in\{0,1\}$.
\endproclaim

\demo{Proof} The first claim follows directly from  Lemma~5.1 and Theorem~5.4 because
the integral in  \thetag{5-1}  equals 
$$
\multline
\sum_{a_1,\dots,a_n=0}^1\sqrt{\kappa(a_1)\lambda(a_1)}\prod_{i=1}^{n-1}\sqrt{Q^{(i)}_{a_{i},a_{i+1}}P^{(i)}_{a_{i},a_{i+1}}}\\
=
\sum_{a_1=0}^1\sqrt{\kappa(a_1)\lambda(a_1)}\prod_{i=1}^{n-1}\sum_{a_{i+1}=0}^1\sqrt{Q^{(i)}_{a_{i},a_{i+1}}P^{(i)}_{a_{i},a_{i+1}}}.
\endmultline
$$
The second claim follows from the fact that $\mu\sim\nu$ and 
$$
\frac{d\nu}{d\mu}(x)=\lim_{n\to\infty}\frac{\nu([x_1,\dots,x_n]_1^n)}{\mu([x_1,\dots,x_n]_1^n)}
$$
 for a.e. $x\in X$.
The final claim was proved in \cite{Lo}.
\qed
\enddemo
We can rewrite \thetag{5-3}  formally in the following form
$$
\sum_{a_1=0}^1\sqrt{\kappa(a_1)\lambda(a_1)}\prod_{i=1}^{\infty}\sum_{a_{i+1}=0}^1\sqrt{Q^{(i)}_{a_{i},a_{i+1}}P^{(i)}_{a_{i},a_{i+1}}}\ne 0,\tag5-4
$$
which is close to the classical Kakutani criterium from \cite{Ka}.

\head 6.  Krieger's type of  tail equivalence~relations equipped with stationary Markov measures
\endhead

In this section we compute Krieger's type of  the tail equivalence relation $\Cal R$ on $X$ equipped  with ``stationary'' Markov measures.

\proclaim{Proposition 6.1} Let $\lambda$ be a non-degenerate distribution on $\{0,1\}$.
Given  a non-degenerate  stochastic matrix  $P$, we let $\boldsymbol P:=(P^{(n)})_{n=1}^\infty$ with $P^{(n)}=P$ for each $n>0$.
Denote by $\mu$ the Markov measure on $X$ determined by $(\lambda,\boldsymbol P)$.
Denote by $\Gamma$ the subgroup of $\Bbb R^*_+$ generated by reals
$\frac{P_{0,0}}{P_{1,1}}$ and $\frac{P_{0,0}^2}{P_{0,1}P_{1,0}}$.
\roster
\item"$\bullet$"
If  $\Gamma=\{1\}$
then $\Cal R$ on $(X,\mu)$ is ergodic and of type $II_1$.
\item"$\bullet$"
If $\Gamma=\{\lambda^n\mid n\in\Bbb Z\}$ for some $\lambda\in(0,1)$
then $\Cal R$ on $(X,\mu)$ is ergodic and of type $III_\lambda$.
\item"$\bullet$"
If $\Gamma$ is dense in $\Bbb R^*_+$
then $\Cal R$ on $(X,\mu)$ is ergodic and  of type $III_1$.
\endroster
\endproclaim
\demo{Proof}
We first note that  $\Cal R$ is ergodic by Theorem~5.4.
The restriction of $\Cal R$ to the subset $[0]_1$ is isomorphic to $\Cal R$.
Hence Krieger's type of $\Cal R\cap ([0]_1\times[0]_1)$ equals Krieger's type of $\Cal R$.
Given $k>0$, we define a transformation $\delta_k$ of $X$ by setting for each $x=(x_j)_{j>0}$ and $j>0$,
$$
(\delta_kx)_j=\cases x_j&\text{if }j\ne k\\
x_k^*&\text{if }j=k,\endcases
$$
where $x_k^*$ is determined by the condition $\{x_k,x_k^*\}=\{0,1\}$.
Then $\delta_k$ is an invertible $\mu$-nonsingular transformation of $X$.
Denote by $\Lambda$ the group of transformations of $X$ generated by $\delta_k$, $k=2,3,\dots$.
Then  two points $x,y\in[0]_1$ are $\Cal R$-equivalent if and only if
$y\in\{\delta x\mid \delta\in\Lambda\}$.
Hence $\Delta_{\Cal R,\mu}(x,y)\in\{\frac{d\mu\circ\delta}{d\mu}(x )\mid \delta\in\Lambda\}$.
It is straightforward to verify that
$$
\frac{d\mu\circ\delta_k}{d\mu}(x )\in\left\{\frac{P_{0,0}}{P_{1,1}},
\frac{P_{1,1}}{P_{0,0}},
\frac{P_{0,0}^2}{P_{0,1}P_{1,0}},
\frac{P_{0,1}P_{1,0}}{P_{0,0}^2},\frac{P_{1,0}P_{0,1}}{P_{1,1}^2},\frac{P_{1,1}^2}{P_{1,0}P_{0,1}}\right\}
$$
for each $x\in X$ and $k>1$.
It follows from this and the cocycle identity that
 the Radon-Nikodym cocycle $\Delta_{\Cal R,\mu}$
restricted to $\Cal R\cap([0]_1\times[0]_1)$ takes  its values in $\Gamma$.
Therefore it suffices to show that $\frac{P_{1,1}}{P_{0,0}}$ and $\frac{P_{0,1}P_{1,0}}{P_{0,0}^2}$ are essential values of this cocycle.
Take a cylinder $[0,a_2,\dots,a_n]_1^n\subset [0]_1$.
Let $N>n$ and let $A:=[0,a_2,\dots,a_n]_1^n\cap[0,0,0]_{N-1}^{N+1}$
and $B:=[0,a_2,\dots,a_n]_1^n\cap[0,0,1]_{N-1}^{N+1}$.
By the Perron-Frobenius theorem, there is a strictly positive vector $(\pi_0,\pi_1)$ such that $\lim_{k\to\infty}(P^{k})_{i,j}=\pi_j$ for each $i,j=0,1$.
Given $\epsilon>0$, if $N$ is large enough then
$$
\align
\frac{\mu(A)}{\mu([0,a_2,\dots,a_n]_1^n)}&=(P^{N-n})_{a_n,0}P_{0,0}^2\ge (1-\epsilon)\pi_0P_{0,0}^2\quad \text{and}\\
\frac{\mu(B)}{\mu([0,a_2,\dots,a_n]_1^n)}&=(P^{N-n})_{a_n,0}P_{0,0}P_{0,1}\ge (1-\epsilon)\pi_0P_{0,0}P_{0,1}.
\endalign
$$
Moreover, $\delta_NA\subset [0,a_2,\dots,a_n]_1^n$,
$\delta_NB\subset [0,a_2,\dots,a_n]_1^n$,
$$
\align
\frac{d\mu\circ\delta_N}{d\mu}(x)&=\frac{P_{0,1}P_{1,0}}{P_{0,0}^2}\quad\text{for all $x\in A$ and }\\
\frac{d\mu\circ\delta_N}{d\mu}(x)&=\frac{P_{0,1}P_{1,1}}{P_{0,0}P_{0,1}}=\frac{P_{1,1}}{P_{0,0}}\quad\text{for all $x\in B$.}
\endalign
$$
Now Lemma~1.1 yields that $\frac{P_{0,1}P_{1,0}}{P_{0,0}^2}$ and
$\frac{P_{1,1}}{P_{0,0}}$ are  essential values of $\Delta_{\Cal R,\mu}\restriction [0]_1$.
\qed
\enddemo

\proclaim{Corollary 6.2} Let $\alpha,\beta\in(0,\infty)$ and let
$
P=\pmatrix
\frac1{1+\alpha} &\frac\alpha{1+\alpha}\\
\frac\beta{1+\beta} & \frac 1{1+\beta}
\endpmatrix.
$
Then $\Cal R$ is of type:
\roster
\item"\rom{(i)}"
 $II_1$ if and only if $\alpha=\beta=1$;
\item"\rom{(ii)}"  $III_\lambda$ for some $\lambda\in(0,1)$ if and only if
there are
integers $p,q\in\Bbb Z$ such that $p\Bbb Z+q\Bbb Z=\Bbb Z$,
$\alpha$ is the positive root of the quadratic equation
$$
\aligned
\alpha^2+(1-\lambda^p)&\alpha-\lambda^q=0\quad\text{and}\\
&\beta=\lambda^{q-p}/\alpha; 
\endaligned
\tag6-1
$$
\item"\rom{(iii)}" 
 $III_1$ if $\Cal R$ is neither of type $II_1$ nor of type $III_\lambda$ for $\lambda\in(0,1)$.
\endroster
\endproclaim
\demo{Proof}
In view of Proposition~6.1,  the claim (i) is obvious.

(ii)
It follows from  Proposition~6.1 that  $\Cal R$ is of type $III_\lambda$ if and only if there are
relatively prime integers $p,q\in\Bbb Z$ such that 
$$
\frac{1+\beta}{1+\alpha}=\lambda^{-p}\quad\text{and}\quad\frac{1+\beta}{(1+\alpha)\alpha\beta}=\lambda^{-q}.
$$
This implies \thetag{6-1}.
We also note that the pair $(\alpha,\beta)=(1,1)$  can not be the solution of \thetag{6-1} because this would imply that  $p=q=0$ and hence $p\Bbb Z+q\Bbb Z\ne\Bbb Z$.

(iii) follows from (i), (ii) and Proposition~6.1.\qed
\enddemo

\head 7. Nonsingular Markov shifts
\endhead

Let  $\mu$ be a non-degenerate and  non-atomic Markov measure  determined by some pair $(\lambda,\boldsymbol P)$.
Let $T$ denote the one-sided shift on $(X,\mu)$.
It is straightforward  to verify that $\mu\circ T^{-1}$ is also a Markov measure on $X$.
This measure is  determined by a pair $(\widehat \lambda,\widehat{\boldsymbol P})$, where
$\widehat \lambda(a)=\sum_{i=0}^1\lambda(i)P^{(1)}_{i,a}$ for $a=0,1$,
and $\widehat{\boldsymbol P}=(\widehat P^{(n)})_{n=1}^\infty$ with $\widehat P^{(n)}=P^{(n+1)}$ for $n>0$.
It follows from Corollary~5.5 that  $T$ is $\mu$-nonsingular  if and only if
$$
\sum_{a_1=0}^1\sqrt{\widehat\lambda(a_1)\lambda(a_1)}\prod_{i=1}^{\infty}\sum_{a_{i+1}=0}^1\sqrt{P^{(i+1)}_{a_{i},a_{i+1}}P^{(i)}_{a_{i},a_{i+1}}}\ne 0,\tag7-1
$$
In this case,  for $\mu$-a.a. $x\in X$, we have
$$
\frac{d\mu\circ T^{-1}}{d\mu}(x)=\frac{\lambda(0)P^{(1)}_{0,x_1}+\lambda(1)P^{(1)}_{1,x_1}}{\lambda(x_1)}\prod_{k=1}^\infty\frac{P^{(k+1)}_{x_{k},x_{k+1}}}{P^{(k)}_{x_k,x_{k+1}}}.\tag7-2
$$

\definition{Definition 7.1}
We call the dynamical system $(X,\mu,T)$ the
{\it nonsingular  one-sided Markov shift} if
 $\mu$ is a non-degenerate and non-atomic Markov measure determined by $(\lambda,\boldsymbol P)$ such that  \thetag{7-1} holds.
\enddefinition

Since $\Cal S_T=\Cal R$, it follows from Theorem~5.4 and Lemma~2.1$(ii)$ that all   nonsingular   one-sided Markov shift are exact (and hence ergodic).
We now describe the natural extensions of  these nonsingular endomorphisms.

\example{Example 7.2} Let $(X,\mu,T)$ be a nonsingular one-sided Markov shift  as above.
Denote by $\widetilde X$, $\widetilde T$ and $\pi$ the same objects  as in
Example~2.4.
Then $\widetilde T$ is the two-sided shift on $\widetilde X$.
To define the corresponding  measure $\widetilde \mu$ on  $\widetilde X$, we first set
$$
Q_{i,j}:=\frac{\lambda(i)P^{(1)}_{i,j}}{\lambda(0)P^{(1)}_{0,j}+\lambda(1)P^{(1)}_{1,j}},\qquad i,j=0,1.\tag{7-3}
$$
Then $Q:=\pmatrix Q_{0,0} & Q_{0,1}\\ Q_{1,0} & Q_{1,1}\endpmatrix$ is a left stochastic matrix, i.e. $\sum_{i=0}^1Q_{i,j}=1$ for each $j\in\{0,1\}$.
We now let
$$
\widetilde\mu([a_{-n},\dots,a_0,\dots,a_n]_{-n}^n)=\bigg(\prod_{-n< i\le 1}Q_{a_{i-1},a_i}\bigg)\lambda(a_1)\prod_{1\le i<n}P^{(i)}_{a_i,a_{i+1}}.
$$
It is straightforward  to verify that $\widetilde T$ is $\widetilde\mu$-nonsingular and
for $\widetilde\mu$-a.e. $\widetilde x=(\widetilde x_k)_{k\in\Bbb Z}\in\widetilde X$, we have
 $$
\frac{d\widetilde\mu\circ\widetilde T}{d\widetilde\mu}(\widetilde x)
=
\frac{Q_{x_1,x_2}\lambda(x_2)}{\lambda(x_1)P^{(1)}_{x_1,x_2}}\cdot
\prod_{k=2}^\infty\frac{ P^{(k-1)}_{x_k,x_{k+1}}}{ P^{(k)}_{x_k,x_{k+1}}}.
\tag7-4
$$
We use here the notation $x:=\pi(\widetilde x)$, $x=(x_n)_{n>0}$ and hence $x_n=\widetilde x_n$ for each $n>0$.
On the other hand, \thetag{7-2} yields that
$$
\frac{d\mu}{d\mu\circ T^{-1}}(Tx)=\frac{\lambda(x_2)}{\lambda(0)P^{(1)}_{0,x_2}+\lambda(1)P^{(1)}_{1,x_2}}\cdot
\prod_{k=2}^\infty \frac{P^{(k-1)}_{x_k,x_{k+1}}}{ P^{(k)}_{x_k,x_{k+1}}}.
$$
This, \thetag{7-3} and \thetag{7-4} yield that $\frac{d\widetilde\mu\circ\widetilde T}{d\widetilde\mu}(\widetilde x)=\frac{d\mu}{d\mu\circ T^{-1}}\circ T(\pi(\widetilde x))$ for  $\widetilde\mu$-a.e. $\widetilde x$, i.e.
$\omega_{\widetilde T,\widetilde\mu}=\omega_{T,\mu}\circ\pi$.
Let $\goth B$ and $\widetilde{\goth B}$ denote the standard Borel $\sigma$-algebras on $X$ and $\widetilde X$ respectively.
Since $\bigvee_{n\in\Bbb Z}\widetilde T^n\pi^{-1}(\goth B)=\widetilde{\goth B}$, it follows that $\widetilde T$ is the natural extension of $T$, as desired.
Since $T$ is exact, $\widetilde T$ is a  $K$-automorphism.
We also deduce from~\thetag{7-4} that for each $n>0$,
$$
\frac{d \widetilde\mu\circ \widetilde T^n}{d\widetilde\mu}(\widetilde x)=
 \frac{\lambda(x_{n+1})}{\lambda(x_1)}\prod_{i=1}^n\frac{Q_{x_i,x_{i+1}}}{P^{(i)}_{x_i,x_{i+1}}}
 \prod_{i=n+1}^\infty
 \frac{P^{(i-n)}_{x_i,x_{i+1}}}{P^{(i)}_{x_i,x_{i+1}}}\tag7-5
$$
at $\widetilde\mu$-a.e. $\widetilde x$.
We will utilize this formula below.
\endexample

\proclaim{Lemma 7.3} Let $(X,\mu,T)$ be a nonsingular one-sided  Markov
shift and let the  measure  $\mu$ be determined by a pair $(\lambda,\boldsymbol P)$.
Then $\alpha_{\omega_{T,\mu}}\restriction\Cal S_T=\Delta_{\Cal S_T,\mu}/\delta$,
where $\delta:\Cal S_T\to\Bbb R^*_+$ is a cocycle of  $\Cal S_T$ given by $\delta(x,y)=\prod_{i=1}^\infty Q_{y_i,y_{i+1}}/Q_{x_i,x_{i+1}}$.
\endproclaim
\demo{Proof}
As in the proof of Theorem~3.1 (see Claim~B) take
 $(x,y)\in S_T$.
  Then there is $n>0$ such that $T^nx=T^ny$, i.e. $x_{i}=y_i$ if $i>n$.
  Let $\widetilde x,\widetilde y\in\widetilde X$ be such that $\pi(\widetilde x)=x$ and
 $\pi(\widetilde y)=y$.
 We have that
 $$
 \alpha_{\omega_{T,\mu}}(x,y)=\frac{ \alpha_{\omega_{T,\mu}}(x,T^nx)}{ \alpha_{\omega_{T,\mu}}(y,T^ny)}=
 \frac{\frac{d \widetilde\mu\circ \widetilde T^n}{d\widetilde\mu}(\widetilde x)}{\frac{d \widetilde\mu\circ \widetilde T^n}{d\widetilde\mu}(\widetilde y)}.
$$
Applying \thetag{7-5} we obtain that
$$
\align
\alpha_{\omega_{T,\mu}}(x,y)
&= \frac{\frac{\lambda(x_{n+1})}{\lambda(x_1)}\prod_{i=1}^n\frac{Q_{x_i,x_{i+1}}}{P^{(i)}_{x_i,x_{i+1}}}
 \prod_{i=n+1}^\infty
 \frac{P^{(i-n)}_{x_i,x_{i+1}}}{P^{(i)}_{x_i,x_{i+1}}}}
 {
 \frac{\lambda(y_{n+1})}{\lambda(y_1)}\prod_{i=1}^n\frac{Q_{y_i,y_{i+1}}}{P^{(i)}_{y_i,y_{i+1}}}
 \prod_{i=n+1}^\infty
 \frac{P^{(i-n)}_{y_i,y_{i+1}}}{P^{(i)}_{y_i,y_{i+1}}}
 }\\
&=
\frac{\lambda(y_1)}{\lambda(x_1)}
\prod_{i=1}^n
\frac{P^{(i)}_{y_i,y_{i+1}}}
{P^{(i)}_{x_i,x_{i+1}}}
 \prod_{i=1}^n
 \frac{Q_{x_i,x_{i+1}}}
  {Q_{y_i,y_{i+1}}}.
\endalign
 $$
Hence
$
 \alpha_{\omega_{T,\mu}}(x,y)=\Delta_{\Cal S_T,\mu}(x,y)/\delta(x,y),
$
as desired.
 \qed
\enddemo

\remark{Remark 7.4} The following assertions are verified straightforwardly.
\roster
\item"$(i)$"
 $\delta$ is trivial if and only if  $Q_{i,j}=0.5$ for all $i,j=0,1$.
This happens if only if
$\lambda(0)=\lambda(1)=0.5$ and
$P^{(1)}_{0,j}=P^{(1)}_{1,j}$ for $j=0,1$.
\item"$(ii)$"
If $P^{(1)}$ is bistochastic then $Q$ is bistochastic if and only if $\lambda(0)=\lambda(1)=0.5$.
In this case we have $Q=P^{(1)}$ and
$$
\delta(x,y)=\left(\frac{P_{0,1}^{(1)}}{P_{0,0}^{(1)}}\right)^{\sum_{i=1}^\infty(|y_{i+1}-y_i|-|x_{i+1}-x_i|)}
$$
for all $(x,y)\in\Cal S_T$.
\item"$(iii)$" $Q_{i,0}=Q_{i,1}$  if and only if $P^{(1)}_{0,i}=P^{(1)}_{1,i}$,  $i=0,1$.
\endroster
\endremark

\proclaim{Theorem 7.5}
Let $(X,\mu,T)$ be a  nonsingular one-sided Markov shift.
If the cocycle $\Delta_{\Cal S_T,\mu}/\delta$ is ergodic then the Maharam extension
of  the natural extension $(\widetilde X,\widetilde\mu,\widetilde T)$ of $T$ is a $K$-automorphism.
If, moreover,  $\widetilde T$ is conservative then $\widetilde T$ is weakly mixing and of type $III_1$.
\endproclaim
\demo{Idea of the proof}
Repeat the argument  in the beginning of Claim~B from the proof of Theorem~3.1 almost literally and then apply Lemma~7.3. \qed
\enddemo

We now prove a necessary  condition for conservativeness of $\widetilde T$.

\proclaim{Lemma  7.6}
Let $P^{(1)}_{0,0}=P^{(1)}_{1,0}$.
If there exist  $\eta>0$ and $k>0$ such that
 $|Q_{b,a}-P^{(n)}_{a,b}|>\eta$  for all $a,b\in\{0,1\}$ and $n>k$ then $\widetilde T$ is not conservative.
 \endproclaim
 \demo{Proof}
By the condition of the lemma and Remark~7.4(iii), there is a probability distribution $q$ on $\{0,1\}$ such that $Q_{a,b}=q(a)$ for all $a,b\in\{0,1\}$.

 It follows from   \thetag{7-5} and Fatou's lemma  that
 $$
 \int_{\widetilde X}  \sqrt{\frac{d \widetilde\mu\circ \widetilde T^n}{d\widetilde\mu}(\widetilde x)}
\,d\widetilde\mu(\widetilde x)
\le
\liminf_{N\to\infty}\int_{\widetilde X}\sqrt{
 \frac{\lambda(x_{n+1})}{\lambda(x_1)}\prod_{i=1}^n\frac{Q_{x_i,x_{i+1}}}{P^{(i)}_{x_i,x_{i+1}}}
 \prod_{i=n+1}^N
 \frac{P^{(i-n)}_{x_i,x_{i+1}}}{P^{(i)}_{x_i,x_{i+1}}}\,}     \,d\widetilde\mu(\widetilde x).
 $$
 The integral in the righthand side of this inequality equals
 $$
 \align
 &\sum_{a_1,\dots,a_{N+1}=0}^1     \sqrt{
\frac{\lambda(a_{n+1})}{\lambda(a_1)}\prod_{i=1}^n\frac{Q_{a_i,a_{i+1}}}{P^{(i)}_{a_i,a_{i+1}}}
 \prod_{i=n+1}^N
 \frac{P^{(i-n)}_{a_i,a_{i+1}}}{P^{(i)}_{a_i,a_{i+1}}}}
\cdot \mu([a_1\dots a_{N+1}]_1^{N+1})\\
&=\sum_{a_1,\dots,a_{N+1}=0}^1\sqrt{
\lambda(a_{n+1}){\lambda(a_1)}\prod_{i=1}^n{Q_{a_i,a_{i+1}}}{P^{(i)}_{a_i,a_{i+1}}}
 \prod_{i=n+1}^N
{P^{(i-n)}_{a_i,a_{i+1}}}{P^{(i)}_{a_i,a_{i+1}}}}\\
&=\sum_{a_1,\dots,a_{n+1}=0}^1\sqrt{
\lambda(a_{n+1}){\lambda(a_1)}\prod_{i=1}^{n}{Q_{a_{i},a_{i+1}}}{P^{(i)}_{a_i,a_{i+1}}}}
 \prod_{i=n+1}^{N}\sum_{a_{i+1}=0}^1
\sqrt{{P^{(i-n)}_{a_i,a_{i+1}}}{P^{(i)}_{a_i,a_{i+1}}}}.
\endalign
$$
Since $\sum_{s=0}^1
\sqrt{{P^{(i-n)}_{u,s}}{P^{(i)}_{u,s}}}\le 1$ and $\lambda(u)<1$ for each $u=0,1$ and $i=n+1,\dots,N$, it follows that
$$
\align
 \int_{\widetilde X}  \sqrt{\frac{d \widetilde\mu\circ \widetilde T^n}{d\widetilde\mu}(\widetilde x)}
\,&d\widetilde\mu(\widetilde x)  
\le
\sum_{a_1,\dots,a_{n+1}=0}^1\sqrt{
\prod_{i=1}^{n}{Q_{a_{i},a_{i+1}}}{P^{(i)}_{a_i,a_{i+1}}}}\\
&=
\sum_{a_1,\dots,a_{n+1}=0}^1\sqrt{q(a_1)}\,\sqrt{
\prod_{i=1}^{n-1}P^{(i)}_{a_i,a_{i+1}}q(a_{i+1})}\,\sqrt{P^{(n)}_{a_n,a_{n+1}}}\\
&\le
2\sum_{a_1,\dots,a_{n}=0}^1\sqrt{
\prod_{i=2}^{n}P^{(i-1)}_{a_{i-1},a_{i}}q(a_{i})}\\
&=2
\sum_{a_1,\dots,a_{k+1}=0}^1\sqrt{
\prod_{i=2}^{k+1}P^{(i-1)}_{a_{i-1},a_{i}}q(a_{i})}\prod_{i=k+2}^n\sum_{a_i=0}^1\sqrt{P^{(i-1)}_{a_{i-1},a_{i}}q({a_{i})}}.
\endalign
$$
Since there is $\xi<1$ such that
$$
\sup\{\sqrt{tq(a)}+\sqrt{(1-t)(1-q(a))}\mid t\in (0,1),\, |t-q(a)|\ge\eta\}\le\xi,
$$
it follows from the condition of the proposition that
$$
\prod_{i=k+2}^n\sum_{a_i=0}^1\sqrt{P^{(i-1)}_{a_{i-1},a_{i}}q({a_{i})}}\le \xi^{n-k-1}.
$$
This yields 
$$
 \int_{\widetilde X}  \sqrt{\frac{d \widetilde\mu\circ \widetilde T^n}{d\widetilde\mu}(\widetilde x)}
\,d\widetilde\mu(\widetilde x) \le 2 \xi^{n-k-1}
\prod_{i=k+2}^n\sum_{a_i=0}^1\sqrt{P^{(i-1)}_{a_{i-1},a_{i}}q(a_i)} \le 2\xi^{n-k-1}.
$$
Therefore $\sum_{n=1}^\infty\langle U_{\widetilde T}^n1,1\rangle<\infty$.
Hence $\widetilde T$ is not conservative by Lemma~2.2.
\qed
\enddemo

\head 8. Maharam extensions of Markov shifts. Bistochastic  case
\endhead

Let $\mu$ be a (non-degenerate and non-atomic) Markov measure determined by a pair $(\lambda,\boldsymbol P)$ for a sequence $\boldsymbol P=(P^{(n)})_{n=1}^\infty$ of bistochastic $2\times 2$ matrices $P^{(n)}$, $n\ge 1$.
It is convenient now to identify $\{0,1\}$ with the group $\Bbb Z/2\Bbb Z$.
Then the space $X$ can be considered as the compact Abelian group $(\Bbb Z/2\Bbb Z)^\Bbb N$.
Let  $\theta:X\to X$ denote the following group  homomorphism 
$$
X\ni x=(x_1,x_2,x_3,\dots)\mapsto(x_1,x_2-x_1,x_3-x_2,\dots) \in X.
$$
This homomorphism has been proved to be useful   in \cite{Do--Qu}.
Of course, $\theta$ is one-to-one and continuous.
For each $y=(y_1,y_2,\dots)\in X$, we have that $\theta(y_1,y_1+y_2,y_1+y_2+y_3,\dots)=y$.
Hence $\theta$ is onto.
Thus $\theta$ is a (topological) automorphism of $X$.
Since $P^{(n)}$ is bistochastic, it follows  that $P^{(n)}_{a+c,b+c}=P^{(n)}_{a,b}$ for all $a,b,c\in\Bbb Z/2\Bbb Z$ and $n>0$.
This yields
$$
\align
\mu\circ\theta^{-1}([y_1,\dots,y_n]_1^n)&=\mu([y_1,y_1+y_2,\dots,y_1+\cdots+y_n]_1^n)\\
&=\lambda(y_1)P^{(1)}_{y_1,y_1+y_2}P^{(2)}_{y_1+y_2,y_1+y_2+y_3}\cdots
P^{(n-1)}_{y_1+\cdots+ y_{n-1},y_1+\cdots+y_n}\\
&=\lambda(y_1)P^{(1)}_{0,y_2}P^{(2)}_{0,y_3}\cdots P^{(n-1)}_{0,y_n}.
\endalign
$$
Hence $\mu\circ\theta^{-1}$ is a Bernoulli measure $\bigotimes_{n=1}^\infty\mu_n$ on $X$, where $\mu_1:=\lambda$ and $\mu_n(i):=P^{(n-1)}_{0,i}$, $i=0,1$ and $n>1$.
We claim that
$$
(\theta\times\theta)(\Cal R)=\Cal R_0:=\left\{(x,y)\in\Cal R\mid \sum_{i=1}^\infty(y_i-x_i)=0\right\}.
\tag8-1
$$
The inclusion $(\theta\times\theta)(\Cal R)\subset\Cal R_0$ is obvious. 
Conversely, if $(x,y)\in\Cal R_0$, then there is $N>0$ such that  $y_1+\cdots +y_N=x_1+\cdots+x_N$ and $x_i=y_i$ for each $i>N$.
Then $(\theta^{-1}x,\theta^{-1}y)=((x_1,x_1+x_2,\dots), (y_1,y_1+y_2,\dots))\in\Cal R$.
Thus \thetag{8-1} is proved.

Let $T$ stand for the one-sided shift on $X$.
It is easy to verify that
$$
\theta T\theta^{-1}(y_1,y_2,y_3,\dots)=(y_1+y_2,y_3,y_4,\dots)
$$
for each $(y_1,y_2,\dots)\in X$. 
Of course, $T$ is $\mu$-nonsingular if and only if $\theta T\theta^{-1}$ is $\mu\circ\theta^{-1}$-nonsingular.
In turn,  the latter holds if and only if $T$ is
$\mu\circ\theta^{-1}$-nonsingular.
Indeed, for $a\in\Bbb Z/2\Bbb Z$, we set $A_a:=\{(u,v)\in (\Bbb Z/2\Bbb Z)^2\mid u+v=a\}$.
Then we have 
$$
\aligned
\frac{d\mu\circ T^{-1}}{d\mu}(\theta^{-1} y)&=\frac{d(\mu\circ\theta^{-1})\circ \theta T^{-1}\theta^{-1}}{d\mu\circ\theta^{-1}}(y)\\
&=\frac{(\mu_1\otimes\mu_2)(A_{y_1})}{\mu_1(y_1)}\prod_{n\ge2}\frac{\mu_{n+1}(y_n)}{\mu_n(y_n)}\\
&=\frac{(\mu_1\otimes\mu_2)(A_{y_1})}{\mu_2(y_1)}\cdot\frac{d(\mu\circ\theta^{-1})\circ T^{-1}}{d(\mu\circ\theta^{-1})}(y)
\endaligned
\tag8-2
$$
for $\mu\circ\theta^{-1}$-a.e. $y=(y_1,y_2,\dots)\in X$.
Therefore  $T$ is $\mu$-nonsingular if and only if $T$ is $\mu\circ\theta^{-1}$-nonsingular.
Equivalently, in the bistochastic case,
 \thetag{7-1} is equivalent to \thetag{2-1}. 
We then call the dynamical system $(X,\mu,T)$  a {\it bistochastic} nonsingular one-sided Markov shift.
Let $(\widetilde X,\widetilde\mu,\widetilde T)$ denote the natural extension of $(X,\mu,T)$.

\proclaim{Theorem 8.1} 
Let $(X,\mu,T)$ be a bistochastic nonsingular one-sided Markov shift as above with $\lambda(0)=\lambda(1)=0.5$ and $P^{(1)}=\pmatrix 0.5 &0.5\\0.5&0.5\endpmatrix$.
Then the following are satisfied.
\roster
\item"$(i)$" $(\widetilde X,\widetilde\mu,\widetilde T)$ is conservative 
if and only if $(\widetilde X,\widetilde{\mu\circ\theta^{-1}},\widetilde T)$ is conservative.
\item"$(ii)$"
The Maharam extension of $(\widetilde X,\widetilde\mu,\widetilde T)$ is a $K$-automorphism 
if 
the cocycle $\Delta_{\Cal R_0,\mu\circ\theta^{-1}}$ is ergodic.
\item"$(iii)$" If $(\widetilde X,\widetilde\mu,\widetilde T)$ is conservative then it is weakly mixing and either of type $II_1$ (if  $\sum_{n\ge 1}(P^{(n)}_{0,0}-0.5)^2<\infty$) or of type $III_1$ (otherwise).
In the latter case,
the Maharam extension of $(\widetilde X,\widetilde{\mu},\widetilde T)$ is a $K$-automorphism.
\endroster
\endproclaim

\demo{Proof} It follows from the assumption of the theorem that $(\mu_1\otimes\mu_2)(A_{a})={\mu_1(a)}$ for each $a\in\Bbb Z/2\Bbb Z$ and $\mu_2=\mu_1$.
Then \thetag{8-2} yields that $\omega_{ T,\mu\circ\theta^{-1}}=\omega_{T,\mu}\circ\theta^{-1}$.
Lemma~2.3(vi) yields now that $T$ is $\mu$-recurrent if and only if $T$ is $\mu\circ \theta^{-1}$-recurrent.
Hence,  $(i)$ follows from Lemma~2.5(i).

It follows from Remark~7.4(i) and the assumption of the theorem that the cocycle $\delta$ is trivial.
Then, by Theorem~7.5, the Maharam extension of $\widetilde T$ is a $K$-automorphism
if the cocycle $\Delta_{\Cal S_T,\mu}$ is ergodic.
It remains to note that  (see \thetag{8-1})
 $$
 \Delta_{\Cal S_T,\mu}(\theta^{-1}x,\theta^{-1}y)=\Delta_{\theta\times\theta(\Cal S_T),\mu\circ\theta^{-1}}(x,y)=\Delta_{\Cal R_0,\mu\circ\theta^{-1}}(x,y)
 $$
 for all $(x,y)\in\Cal S_T$.
 Thus $(ii)$ is proved.
 
Since $\Cal R_0$ contains the equivalence relation generated by the group of  finite permutations of $\Bbb N$ acting on $X$ in the natural way, it follows from
 Remark 1.7 that   $\Delta_{\Cal R_0,\mu\circ\theta^{-1}}$ is ergodic 
if and only if  $\Delta_{\Cal R,\mu\circ\theta^{-1}}$ is ergodic.
Therefore arguing as in the proof of Theorem~3.1 (see also Corollary~3.3) and utilizing
 $(i)$ and $(ii)$ we prove~$(iii)$.
\qed
\enddemo

We note that, under the condition of the theorem, if  $(\widetilde X,\widetilde\mu,\widetilde T)$ is conservative and of type $II_1$ then $\mu$ is $T$-cohomologous to a Markov measure determined by the pair $(\lambda, (P^{(1)})_{n=1}^\infty)$.
Thus   $\mu$ is $T$-cohomologous to  a Bernoulli (i.e. infinite product) measure.
The converse follows from the proposition below:

\proclaim{Proposition 8.2}
Let $(X,\mu,T)$ be as in Theorem~8.1.
If $\mu$ is equivalent to a Bernoulli measure
 then $\widetilde T$ is ergodic and of type $II_1$.
\endproclaim

\demo{Proof}
Let $\mu$ be equivalent to a Markov measure determined by a pair $(\rho, (V^{(n)})_{n=1}^\infty)$ with $V^{(n)}_{0,0}=V^{(n)}_{1,0}$ for all $n>0$.
Then $V^{(n)}$ is close to $P^{(n)}$ for all sufficiently large $n$ by Corollary~5.5.
Since $P^{(n)}$ is bistochastic, it follows that $V^{(n)}_{0,0}\to 0.5$ as $n\to \infty$.
If follows that there is $\delta>0$ such that $\inf_{n>0}\min_{i,j}V^{(n)}_{i,j}>\delta$.
Hence there is $\delta'>0$ such that $\inf_{n>0}\min_{i,j}P^{(n)}_{i,j}>\delta'$.
It follows from Corollary~5.5 that $\sum_{n=1}^\infty(P^{(n)}_{0,0}-V^{(n)}_{0,0})^2<\infty$ 
and $\sum_{n=1}^\infty(1-P^{(n)}_{0,0}-V^{(n)}_{0,0})^2<\infty$.
Therefore $\sum_{n=1}^\infty(P^{(n)}_{0,0}-0.5)^2<\infty$ and hence
$\sum_{n=1}^\infty(P^{(n)}_{i,j}-0.5)^2<\infty$
for all $i,j\in\{0,1\}$.
It remains to apply Corollary~5.5 again (but in the opposite direction).
\qed
\enddemo

\remark{Remark 8.3} Utilizing Theorem 8.1 we can produce concrete examples of type $III_1$ ergodic conservative natural extensions of bistochastic nonsingular  one-sided Markov shifts.
For that take a concrete example  of a nonsingular one-sided Bernoulli shift $(X,\bigotimes_{n=1}^\infty\mu_n, T)$ whose natural extension is conservative ergodic and of type $III_1$ 
such that $\mu_1(0)=\mu_2(0)=0.5$ (such systems were constructed in Theorem~4.2 and earlier in \cite{Kr}, \cite{Ha}, \cite{Ko1} and \cite{Ko2}) and set $\mu:=(\bigotimes_{n=1}^\infty\mu_n)\circ\theta$.
Then the system $(\widetilde X,\widetilde\mu,\widetilde T)$ is as desired.
We note also that $\mu$ is not equivalent to any Bernoulli measure on $X$ according to Proposition~8.2.
\endremark

\head 9. Maharam extensions of  Markov shifts. General case
\endhead

Let $J$ be an infinite  subset of $\Bbb N$ such that its complement $\Bbb N\setminus J$ is also infinite.
We endow $J$  and $\Bbb N\setminus J$  with the induced (from $\Bbb N$) linear ordering.
Using this ordering we may identify naturally $\{0,1\}^J$ with $X$
and $\{0,1\}^{\Bbb N\setminus J}$
 with $X$.
We denote by $\phi_J:X\ni x\mapsto x|J\in X$  and
$\phi_{\Bbb N\setminus J}:X\ni x\mapsto x|(\Bbb N\setminus J)\in X$ the corresponding restriction maps.
Then the mapping
$$
X\ni x\mapsto (\phi_J\times\phi_{\Bbb N\setminus J})(x)\in X\times X\tag9-1
$$
is a homeomorphism.
Moreover, the mapping $(\phi_J\times\phi_{\Bbb N\setminus J})\times (\phi_J\times\phi_{\Bbb N\setminus J})$ maps bijectively $\Cal R$
onto $\Cal R\times\Cal R$.

Let $\mu$ be a probability measure on $X$ such that $\Cal R$ is $\mu$-nonsingular.
In view of the identification \thetag{9-1}, we may consider $\mu$
as a measure on $X\times X$. Then $\Cal R\times\Cal R$ is $\mu$-nonsingular.
Now
we let $\mu_J:=\mu\circ\phi_J^{-1}$.
Then $\mu$ admits a disintegration
 $$
\mu=\int_X\delta_x\times\mu^{(x)}\,d\mu_J(x)\tag9-2
$$
  relative to $\mu_J$, where $\delta_x$ is the Kronecker measure supported at $x$ and $X\ni x\mapsto\mu^{(x)}$ is the  corresponding canonical system of conditional measures on $X$.
  Since $\Cal R\times\Cal R$ is $\mu$-nonsingular, $\Cal R$ is $\mu^{(x)}$-nonsingular for $\mu_J$-a.a. $x\in X$.
  Moreover, $\mu^{(x)}$ and $\mu^{(x')}$ are equivalent whenever $(x,x')\in\Cal R$
   and
  $$
  \Delta_{\Cal R\times\Cal R,\mu}((x,y),(x',y'))=\Delta_{\Cal R,\mu_J}(x,x')\frac{d\mu^{(x')}}{d\mu^{(x)}}(y)\Delta_{\Cal R,\mu^{(x')}}(y,y').
\tag{9-3}
  $$

\proclaim{Lemma 9.1}
Suppose that $\Cal R\times\Cal R$ is $\mu$-ergodic.
If $(\Cal R,\mu^{(x)})$ is ergodic and of type $III_1$ for $\mu_J$-a.e. $x\in X$ then $(\Cal R\times\Cal R,\mu)$ is also of type $III_1$.
\endproclaim
\demo{Proof}
Denote the skew product equivalence relation $(\Cal R\times\Cal R)(\Delta_{\Cal R\times\Cal R,\mu})$ on the product space $(X\times X\times \Bbb R^*_+,\mu\times\lambda_{ \Bbb R^*_+})$ by
$\widetilde{\Cal R}$.
We need to show that $\widetilde {\Cal R}$ is ergodic.
Let $A$ be an $\widetilde {\Cal R}$-invariant subset of $X\times X\times \Bbb R^*_+$.
Given $x\in X$, we let
$$
A_x:=\{(y,z)\in X\times \Bbb R^*_+\mid (x,y,z)\in A\}.
$$
Since $A$ is  $\widetilde {\Cal R}$-invariant, it follows from \thetag{9-3} that $A_x$ is $\Cal R(\Delta_{\Cal R,\mu_x})$-invariant
 for $\mu_J$-a.e. $x\in X$.
Indeed,  $(x,y)\sim(x,y')$ whenever $(y,y')\in\Cal R$
and hence \thetag{9-3} yields
$$
 \Delta_{\Cal R\times\Cal R,\mu}((x,y),(x,y'))=\Delta_{\Cal R,\mu^{(x)}}(y,y').
$$
 By the assumptions of the lemma,  $(\Cal R,\mu^{(x)})$ is ergodic and of type $III_1$.
Hence the skew product extension
 $\Cal R(\Delta_{\Cal R,\mu^{(x)}})$ is ergodic.
Therefore we have either $(\mu_x\times\lambda_{ \Bbb R^*_+})(A_x)=0$ or
 $(\mu_x\times\lambda_{ \Bbb R^*_+})((X\times\Bbb R^*_+)
\setminus A_x)=0$ for a.e. $x\in X$.
 We let
$$
\align
B&:=\{x\in X\mid (\mu_x\times\lambda_{ \Bbb R^*_+})((X\times\Bbb R^*_+)\setminus A_x)=0\}\quad
\text{
 and}\\
 \widetilde B&:=\{(x,y,z)\in X\times X\times \Bbb R^*_+ \mid x\in B\}.
\endalign
$$
 By the Fubini theorem, $(\mu\times\lambda_{ \Bbb R^*_+})(A\triangle \widetilde B)=0$.
 Hence $\widetilde B$ is $\widetilde{ \Cal R}$-invariant (mod $0$).
 This is possible if and only if $B$ is $\Cal R$-invariant.
Since $\Cal R\times\Cal R$ is $\mu$-ergodic by the assumptions of the lemma,
it follows that $\Cal R$ is $\mu_J$-ergodic.
This yields   that  either $\mu_J(B)=0$ or
 $\mu_J(X\setminus B)=0$.
 Hence $A$ is either $(\mu\times\lambda_{ \Bbb R^*_+})$-null or  $(\mu\times\lambda_{ \Bbb R^*_+})$-conull, as desired.
 \qed
 \enddemo

In a similar way one can prove a  general statement on ergodicity of cocycles on $\Cal R\times\Cal R$.

\proclaim{Lemma 9.2} Let $G$ be a locally compact second countable group.
Let $\mu$ be a probability measure on $X\times X$ and let $\mu=\int_X\delta_x\times\mu^{(x)}d\kappa(x)$ be a disintegration of $\mu$ with respect to the projection $X\times X\to X$ onto the first coordinate\footnote{Thus the measure $\kappa$ is the projection of $\mu$ to the first coordinate and $X\ni x\mapsto\mu^{(x)}$ is the corresponding system of conditional measures on $X$.}.
Suppose that $\Cal R\times\Cal R$ is $\mu$-ergodic.
Given a cocycle $\alpha:\Cal R\times\Cal R\to G$, we define a measurable field of cocycles  $\alpha^{(x)}:\Cal R\to G$ defined on $(X,\mu^{(x)},\Cal R)$, $x\in X$, by setting $\alpha^{(x)}(y,y'):=\alpha((x,y),(x,y'))$ for all $(y,y')\in\Cal R$.
If  $\alpha^{(x)}$ is ergodic  for $\kappa$-a.e. $x\in X$ then $\alpha$ is also  ergodic.
\endproclaim

We now compute $\mu_J$ and  the conditional measures $\mu^{(x)}$  in the disintegration~\thetag{9-2} in the case  where $\mu$ is a Markov measure and $J=\{1,3,5,\dots\}$.
We also describe the measurable field of cocycles $\delta^{(x)}$, $x\in X$, corresponding to the cocycle $\delta$ defined in Lemma~7.3.

  \proclaim{Lemma 9.3}
Let  $\mu$ be a Markov measure determined by a pair $(\lambda,\boldsymbol P)$.
Let $J$ be the set of odd positive integers.
Set $\kappa_{u,v}:=\frac{Q_{u,1}Q_{1,v}}{Q_{u,0}Q_{0,v}}$ for $u,v=0,1$.
 Then 
 \roster
 \item"$(i)$" $\mu_J$ is the Markov measure on $X$ determined by the pair $(\lambda,\boldsymbol P_J)$, where $\boldsymbol P_J=(P^{(2n-1,2n)})_{n=1}^\infty$,
 \item"$(ii)$"
  $\mu^{(x)}$ is the Bernoulli measure $\bigotimes_{n=1}^\infty\mu^{(x)}_n$ on $X$, where
 $$
\mu^{(x)}_n(i)=\frac{P^{(2n-1)}_{x_n,i}P^{(2n)}_{i,x_{n+1}}}{P^{(2n-1,2n)}_{x_n,x_{n+1}}}, \quad i=0,1,
\tag9-4
$$
and
 \item"$(iii)$"  $\delta^{(x)}(a,b)=\prod_{i\ge 1}(\kappa_{x_i,x_{i+1}})^{b_i-a_i}$ for $(a,b)\in\Cal R$, $a=(a_i)_{i=1}^\infty$, $b=(b_i)_{i=1}^\infty$,  for $\mu_J$-a.e. $x=(x_1,x_2,\dots)\in X$.
\endroster
  \endproclaim

  \demo{Proof}
  The first claim is obvious.
  To prove the second one, we take
    $n>0$ and a cylinder $[a_1,\dots,a_n]_1^n$.
    Then we have
  $$
  \align
  \mu^{(x)}([a_1,\dots,a_n]_1^n)
  &=\frac{\mu([x_{1},a_1,x_2,a_2,\dots,x_n,a_{n},x_{n+1}]_1^{2n+1})}{\mu([x_1]_{1}\cap[x_2]_3\cap\dots\cap[x_{n+1}]_{2n+1})}
  \\
  &=\frac{\lambda(x_1)P^{(1)}_{x_1,a_1}P^{(2)}_{a_1,x_2}\cdots
P^{(2n-1)}_{x_n,a_n} P^{(2n)}_{a_n,x_{n+1}}}
  {\lambda(x_1)P^{(1,2)}_{x_1,x_2}\cdots P^{(2n-1,2n)}_{x_{n},x_{n+1}}}
\\
&=\frac{P^{(1)}_{x_1,a_1}P^{(2)}_{a_1,x_2}}{P^{(1,2)}_{x_1,x_2}}\cdots
\frac{P^{(2n-1)}_{x_n,a_n} P^{(2n)}_{a_n,x_{n+1}}   }{P^{(2n-1,2n)}_{x_n,x_{n+1}}},
  \endalign
  $$
as desired.
It is straightforward to verify the third
claim.
\qed
\enddemo

Given  $x\in X$ and  $n>0$, let
$$
r^{(x)}_n:=\frac{\mu^{(x)}_n(1)}{\mu^{(x)}_n(0)}=\frac{P^{(2n-1)}_{x_n,1}P^{(2n)}_{1,x_{n+1}}}{P^{(2n-1)}_{x_n,0}P^{(2n)}_{0,x_{n+1}}}.
$$
Denote by $\Cal L^{(x)}$ the set of limit points of the sequence $(r_n^{(x)})_{n=1}^\infty$.
Of course, if $(x,y)\in\Cal R$ then $\Cal L^{(x)}=\Cal L^{(y)}$.
Therefore the map $X\ni x\mapsto\Cal L^{(x)}$ is a Borel $\Cal R$-invariant map from $X$ to the space of closed subsets of the ray $[0,+\infty)$ if we endow this space with the Fell topology \cite{Fel}.
Since $\Cal R$ is $\mu_J$-ergodic, there exists a closed subset $\Cal L\subset[0,+\infty)$
such that $\Cal L^{(x)}=\Cal L$  for $\mu_J$-a.e. $x$.
We say that a point $\alpha\in\Cal L^{(x)}$ is {\it good} if there is an increasing sequence $n_1<n_2<\cdots$ of positive integers such that  $\lim_{j\to\infty}r^{(x)}_{n_j}=\alpha$ and $\sum_{i>0}(\alpha-r^{(x)}_{n_j})^2=+\infty$.
If there is a segment  $[\beta,\gamma]\subset[0,+\infty)$ such that  the intersection $\Cal L\cap[\beta,\gamma]$ is infinite then $\Cal L$ contains a  good point.

\proclaim{Theorem 9.4} 
 Let $\mu$ be a Markov measure determined by a pair $(\lambda,\boldsymbol P)$ and 
 let $(X,\mu,T)$ stand for the corresponding  nonsingular one-sided Markov shift.
 \roster
 \item"$(i)$"
 If $\sum_{u=0}^1\lambda(u)P^{(1)}_{u,v}=\lambda(v)$ for $v=0,1$ and 
$\sum_{u=0}^1\sum_{n=1}^\infty(P^{(n)}_{u,0}-P^{(1)}_{u,0})^2<\infty$ then $(\widetilde X,\widetilde\mu,\widetilde T)$ is isomorphic to the measure preserving two-sided  shift on the probability space $(\widetilde X,\widetilde\nu)$, where $\widetilde\nu$ is the Markov measure determined by the pair $(\lambda, (P^{(1)})_{n=1}^\infty)$.
 \item"$(ii)$"
 If  the natural extension $(\widetilde X,\widetilde\mu,\widetilde T)$ of $T$ is conservative
and $\Cal L$ contains a good non-zero point then
the Maharam extension of $(\widetilde X,\widetilde\mu,\widetilde T)$ is a weakly mixing $K$-automorphism.
In particular, $(\widetilde X,\widetilde\mu,\widetilde T)$ is of type $III_1$.
\endroster
\endproclaim

\demo{Proof} (i) follows form the final claim of Corollary~5.5.

(ii)
To show that the Maharam extension of $(\widetilde X,\widetilde\mu,\widetilde T)$ is a $K$-automorphism we will apply Theorem~7.5.
For that we have to prove that the cocycle $\Delta_{S_T,\mu}/\delta$ is ergodic (the cocycle $\delta$ is defined in Lemma~7.3).
Let $J$ stand for the set of odd positive integers.
By Lemma~9.2,  $\Delta_{S_T,\mu}/\delta$ is ergodic if
the cocycle
$
\Delta_{\Cal R,\mu^{(x)}}/\delta^{(x)}
$
of 
 $\Cal R$ on $(X,\mu^{(x)})$ is ergodic  for $\mu_J$-a.e. $x$.
 Given $x\in X$, partition $\Bbb N$ into four subsets $\Bbb N=J^{(x)}_{0,0}\sqcup J^{(x)}_{0,1}\sqcup J^{(x)}_{1,0}\sqcup J^{(x)}_{1,0}$, where $j\in J^{(x)}_{u,v}$ if $x_j=u$ and $x_{j+1}=v$. 
It follows from Lemma~9.3 that the  quadruple $(X,\mu^{(x)},\Cal R,\delta^{(x)})$ splits into direct product
$$
\bigotimes_{u,v=0}^1\left(\{0,1\}^{J^{(x)}_{u,v}},\mu^{(x)}_{u,v},\Cal R_{u,v}^{(x)},\delta^{(x)}_{u,v}\right),
$$
where $\mu^{(x)}_{u,v}:=\bigotimes_{i\in J^{(x)}_{u,v}}\mu^{(x)}_i$, $\Cal R_{u,v}^{(x)}$ is the tail equivalence relation on the product space $\{0,1\}^{J^{(x)}_{u,v}}$
and $\delta^{(x)}_{u,v}$ is a cocycle of $\Cal R_{u,v}^{(x)}$ given by
$$
\delta^{(x)}_{u,v}(a,b):=
(\kappa_{u,v})^{\sum_{i\in J^{(x)}_{u,v}}{(b_i-a_i)}}
$$
Hence it suffices to show that  for some pair $(u,v)$, the cocycle $\Delta_{\Cal R_{u,v}^{(x)},\mu^{(x)}_{u,v}}/\delta_{u,v}^{(x)}$
is ergodic.
It follows from the condition of the theorem that there is
$\alpha\in\Cal L^{(x)}$  and  an increasing sequence $n_1<n_2<\cdots$ of positive integers such that  $\alpha>0$, $\lim_{j\to\infty}r^{(x)}_{n_j}=\alpha$ and $\sum_{i>0}(\alpha-r^{(x)}_{n_j})^2=+\infty$.
Passing to a subsequence, if necessary, we may assume without loss of generality that there is a pair $(u,v)\in\{0,1\}^2$ such that $n_j\in J^{(x)}_{u,v}$ for all $j>0$.
Then $\Delta_{\Cal R_{u,v}^{(x)},\mu^{(x)}_{u,v}}/\delta_{u,v}^{(x)}$
is ergodic in view of~Remark~1.6.
\qed
\enddemo

\head  10. Concluding remarks and open problems
\endhead

\roster

\item Let $\widetilde T$ be the natural extension of a nonsingular one-sided Markov shift.
Suppose that $\widetilde T$ is conservative and not of type $II_1$.
Is the Maharam extension of  $\widetilde T$  a $K$-automorphism?
Theorems~3.1(ii), 8.1(iii) and 9.4(ii) provide only partial answers to this question.

\item  Let $T$ and $R$ be two nonsingular one-sided Bernoulli shifts on the infinite product spaces $(\{0,1\}^\Bbb N,\bigotimes_{n\ge 1}\mu_n)$ and $(\{0,1\}^\Bbb N,\bigotimes_{n\ge 1}\nu_n)$, respectively.
Suppose that the natural extensions $\widetilde T$ and $\widetilde R$ of $T$ and $R$, respectively, are conservative.
Is it possible that $\mu_1(0)\not\in\{\nu_1(0),\nu_1(1)\}$ but $\widetilde T$ and $\widetilde R$ are conjugate as nonsingular transformations?
In particular, can the natural extension of an equilibrial   one-sided Bernoulli shift be conjugate with  the natural extension of a  non-equilibrial  Bernoulli shift? 
We note that in the probability preserving case, i.e. in the case where $\mu_1=\mu_n$ and $\nu_1=\nu_n$ for all $n>1$,  if $\mu_1(0)\not\in\{\nu_1(0),\nu_1(1)\}$ then $h(\widetilde T)\ne h(\widetilde R)$ and hence $\widetilde T$ and $\widetilde R$ are not conjugate.

\item Let $T$ be  a nonsingular one-sided Markov shift on $\{0,1\}^\Bbb N$ and the corresponding Markov measure on this space is not equivalent to a Bernoulli measure.
Suppose that  the natural extension $\widetilde T$ of $T$ is conservative.
Then $\widetilde T$ is weakly mixing.
Are there nonsingular one-sided Bernoulli shifts  whose natural extensions are conjugate to $\widetilde T$?

\item What are the critical dimensions  (see \cite{DanSi} for the definition) of the natural extensions of non-equilibrial nonsingular Bernoulli shifts and  bistochastic Markov shifts?
Can we distinguished  between equilibrial and non-equilibrial Bernoulli shifts using critical dimensions?
Some estimations for these invariants  were obtained in \cite{DoMor} for the equilibrial Bernoulli case. 

\item 
Are there nice criteria for conservativeness  of the natural extensions of the one-sided nonsingular Bernoulli and Markov shifts?

\item
We recall that given an invertible nonsingular transformation $T$ of the standard non-atomic $\sigma$-finite measure space $(X,\mu)$, the {\it ergodic index} $e(T)$ of $T$ is the smallest positive integer $d$ such that the $d$-th Cartesian power $T^{\otimes d}$ of $T$ is not ergodic.
If no such integer exists, $T$ is said to be of infinite ergodic index.
In a similar way one can define {\it index of conservativeness} $c(T)$ for $T$.
Of course,  $e(T)\le c(T)$ (see a survey \cite{DanSi} for more information about these indices).
Suppose now that $ T$ is  
the natural extension of a nonsingular one-sided Bernoulli shift $S$ on a product space $(Y,\nu)$.
It is easy to verify that for each $d>0$, the dynamical system $(X^d,\mu^{\otimes d}, T^{\otimes d})$ is the natural extension of $(Y^d,\nu^{\otimes d}, S^{\otimes d})$.
It is easy to see that $S^{\otimes d}$ is also exact.
Hence  $S^{\otimes d}$ is ergodic.
Therefore, if $T^{\otimes d}$ is conservative then it is ergodic by Lemma~2.5(iii) and (i).
It follows that $e(T)=c(T)$.
A question arises: what are possible values of $e(T)$ when $T$ runs over the set of natural extensions of  all conservative non-singular one-sided Bernoulli  shifts of type $III_1$? 

\item  Let $T$ be an ergodic conservative invertible transformation of type $III_1$.
Then the Maharam extension $\widehat T$ of $T$ is an ergodic conservative transformation
of type $II_\infty$.
Of course, $e(\widehat T)\le e(T)$ and $c(\widehat T)\le c(T)$.
What are possible values of $e(\widehat T)$  and $c(\widehat T)$ when $T$ runs  over the set of natural extensions of  all conservative non-singular one-sided Bernoulli (or Markov) shifts? 
This is related to Z.~Kosloff's question: does the condition $c( T)>2$ imply that
$c(\widehat T)>2$?

\item How ``huge''  is the centralizer\footnote{The {\it centralizer} of an invertible nonsingular transformation $S$ of a probability space is the group of all nonsingular transformations  of this space that commute with $S$.} of the type $III_1$ natural extension $\widetilde T$ of a nonsingular one-sided Bernoulli shift $T$? Does the second centralizer\footnote{The {\it second centralizer} of an invertible nonsingular transformation $S$ of a probability space  is the group of nonsingular  transformations of this space that commute with every nonsingular transformation commuting with $S$.} of $\widetilde T$ consist of just the powers of $\widetilde T$ as in the finite measure preserving case (cf. \cite{Ru})?

\item There are many works devoted to computation  of Krieger's type and the associated flow of the tail equivalence relation  $\Cal R$ on the infinite product of finite spaces equipped with product measures (i.e. odometers of product type). We refer to \cite{DanSi} for references and definition of the associated flow.
It would be interesting to consider similar problems in the Markov case, for instance, in the simplest case, where $X=\{0,1\}^\Bbb N$ and $\mu$ is a nonatomic non-degenerate Markov measure determined by some pair $(\lambda,\boldsymbol P)$.
In particular, is it true that the associated flow of $(X,\mu,\Cal R)$ is AT (i.e. {\it approximatively transitive} in the sense of Connes and Woods  \cite{CoWo})?
Some results  were obtained in \cite{Do--Qu} for the bistochastic case and in Section~6 of the present paper for the simplest  ``stationary'' case.

\item We note\footnote{
The authors thanks Z.~Kosloff for this remark.}
 that the conservative infinite measure preserving $K$-automorp\-hisms that are Maharam extensions of (the natural extensions of) type $III_1$ Bernoulli or Markov shifts do not belong to the class of conservative infinite measure preserving $K$-automorphisms introduced by Parry in \cite{Pa}.
Parry's class consists of irreducible aperiodic recurrent Markov chains with stationary transition probabilities and {\it countable} state space.
Indeed, Aaronson showed in \cite{Aa} that those Markov chains are non-squashable, i.e. every element from the centralizer of such a map $T$ preserves the $T$-invariant measure. 
In contrast, every Maharam extension is squashable. 
\endroster

\enddocument